\documentclass[preprint,5p,twocolumn]{elsarticle}

\usepackage{amssymb,amsmath,amsthm}
\usepackage{bbold} 
\usepackage{xcolor}
\usepackage{hyperref}

\usepackage{pgfplots}
\pgfplotsset{compat=newest}
\usepgfplotslibrary{patchplots,colorbrewer}

\usepackage{pgfmath,pgfplotstable,pgffor}

\usepackage{tikz}

\usepackage{booktabs}
\usepackage{pifont}

\newcommand{\tr}{\intercal}
\newcommand{\one}{\mathbf{1}}
\newcommand{\nf}{\mathrm{f}}
\newcommand{\nE}{\mathrm{e}}
\newcommand{\nv}{\mathrm{v}}
\newcommand{\cvx}[1]{\mathrm{conv}\left(#1\right)}
\newcommand{\cone}[1]{\mathrm{cone}\left(#1\right)}
\newcommand{\interior}[1]{\mathbf{int}\left(#1\right)}
\newcommand{\defeq}{:=}
\newcommand{\order}[1]{\mathbf{O}\left(#1\right)}
\newcommand{\support}[1]{s\left[#1\right]}

\newcommand{\Zono}{\mathcal Z}
\newcommand{\Poly}{\mathcal P}
\newcommand{\Vol}{\mathrm{Vol}}
\newcommand{\Lvol}{\Phi}
\newcommand{\ZLvol}{\Theta}

\newcommand{\SB}[5]{\node at (#1,#2) {\footnotesize \textbf{Section}~\ref{#4}};
\node at (#1,#3) {\footnotesize #5};}

\newcommand{\SBT}[7]{\node at (#1,#2) {\footnotesize \textbf{Section}~\ref{#4}};
\node at (#1,#3) {\footnotesize #5};
\node at (#1,#6) {\footnotesize #7};
}

\begin{document}

\begin{frontmatter}


\title{Polyhedral Control Design: Theory and Methods}

\author[1,3]{Boris Houska}
\ead{borish@shanghaiTech.edu.cn}
\author[2]{Matthias A.~M\"uller}
\ead{mueller@irt.uni-hannover.de}
\author[3]{Mario Eduardo Villanueva\corref{cor1}}
\ead{me.villanueva@imtlucca.it}
\cortext[cor1]{Corresponding author}
\affiliation[1]{organization={ShanghaiTech University},
country={China}}
\affiliation[2]{organization={Leibniz University Hannover},
country={Germany}}
\affiliation[3]{organization={IMT School for Advanced Studies Lucca},
country={Italy}}

\begin{abstract}
In this article, we survey the primary research on polyhedral computing methods for constrained linear control systems. Our focus is on the modeling power of convex optimization, featured to design set-based robust and optimal controllers. In detail, we review the state-of-the-art techniques for computing geometric structures such as robust control invariant polytopes. Moreover, we survey recent methods for constructing control Lyapunov functions with polyhedral epigraphs as well as the extensive literature on robust model predictive control. The article concludes with a discussion of both the complexity and potential of polyhedral computing methods that rely on large-scale convex optimization.
\end{abstract}


\end{frontmatter}

\section{Introduction}
\label{sec::introduction}
Modern robust control theory and methods routinely rely on tools from the field of set-based computing. This refers to both the computation of safety margins for learning in uncertain environments as well as robust feedback control synthesis. A key challenge in this domain, however, is developing tools that can effectively represent and optimize potentially high-dimensional sets. This is crucial in the context of modern, large-scale control synthesis, particularly in the context of data-driven optimal control, where significant uncertainty is present.

In this survey, we argue that polyhedral and polytopic computing tools are among the most promising candidates for representing and optimizing high-dimensional convex sets, thereby addressing the demands of modern set-based control. This assertion is rooted in the fact that polyhedral computing methods leverage linear programming and related convex optimization techniques, which have proven to be versatile and highly scalable. Consequently, these tools can be used to construct and optimize polyhedra and polytopes with thousands or even millions of facets and vertices, making them well-suited for representing and optimizing convex sets.

It is worth noting that polyhedra and polytopes have a rich history, with mathematical theories largely developed over the past two centuries, though early explorations, like Pythagoras' work on triangles, date back to antiquity. Their study stems from diverse motivations, encompassing mathematical curiosity, applications in the complexity analysis of simplex algorithms, computational geometry and topology, global and combinatorial optimization, convex relaxation, set-based verification, and set arithmetics, as well as a long list of downstream tasks in the applied sciences, including robust and optimal control. Given this broad scope and cross-disciplinary development, comprehending the overarching landscape can be challenging. Therefore, this article aims to first survey key milestones in polyhedral computing and computational geometry, followed by examining how these methods facilitate formulating convex optimization problems for robust control tasks.

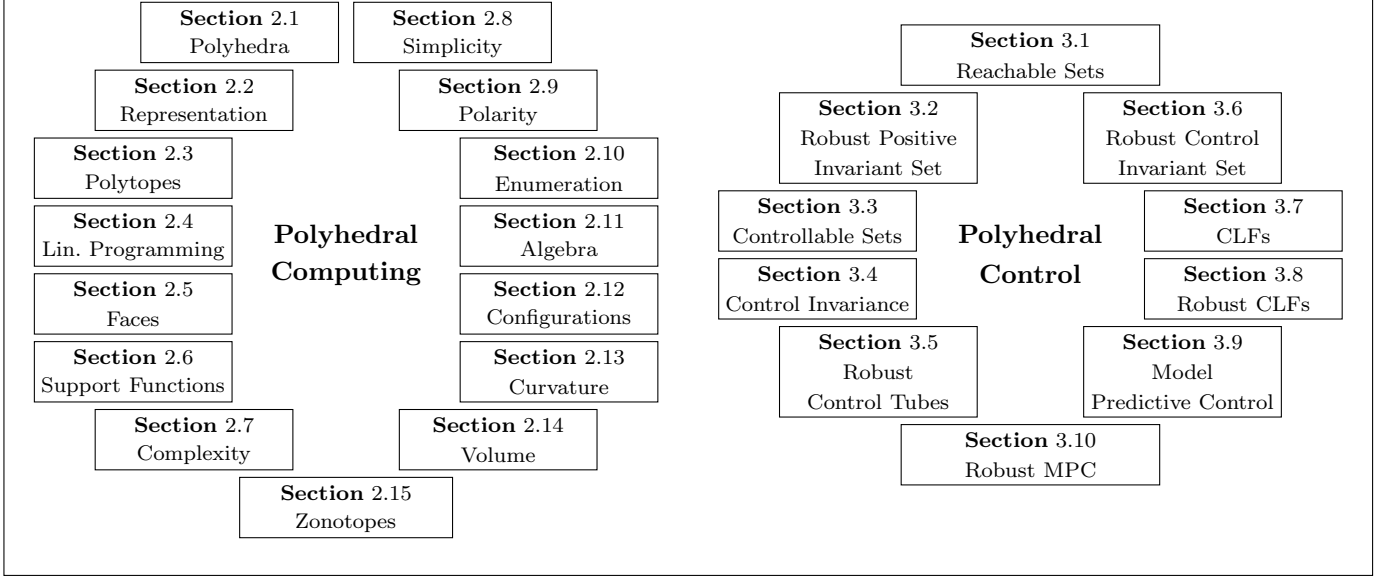
\begin{figure*}
\begin{center}
\begin{tikzpicture}
\draw (0,0) rectangle (18,9);

\node at (9,8.5) {\textbf{Polytopic Computing Theory: The Pathway to Set Based Control via Convex Optimization}};

\draw (1.8,6.8) rectangle (4.4,7.6);
\SB{3.1}{7.4}{7.0}{sec::polyhedra}{Polyhedra}
\draw (4.6,6.8) rectangle (7.2,7.6);
\SB{5.9}{7.4}{7.0}{sec::simplicity}{Simplicity}
\draw (1.2,5.9) rectangle (3.8,6.7);
\SB{2.5}{6.5}{6.1}{sec::representation}{Representation}
\draw (5.2,5.9) rectangle (7.8,6.7);
\SB{6.5}{6.5}{6.1}{sec::polarity}{Polarity}
\draw (0.4,5.0) rectangle (3.0,5.8);
\SB{1.7}{5.6}{5.2}{sec::polytope}{Polytopes}
\draw (6.0,5.0) rectangle (8.6,5.8);
\SB{7.3}{5.6}{5.2}{sec::enumeration}{Enumeration}
\draw (0.4,4.1) rectangle (3.0,4.9);
\SB{1.7}{4.7}{4.3}{sec::LP}{Lin.~Programming}


\node at (4.5,4.5) {\textbf{Polyhedral}};
\node at (4.5,4.0) {\textbf{Computing}};

\node at (13.5,4.5) {\textbf{Polyhedral}};
\node at (13.5,4.0) {\textbf{Control}};

\draw (6.0,4.1) rectangle (8.6,4.9);
\SB{7.3}{4.7}{4.3}{sec::algebra}{Algebra}
\draw (0.4,3.2) rectangle (3.0,4.0);
\SB{1.7}{3.8}{3.4}{sec::faces}{Faces}
\draw (6.0,3.2) rectangle (8.6,4.0);
\SB{7.3}{3.8}{3.4}{sec::configurations}{Configurations}
\draw (0.4,2.3) rectangle (3.0,3.1);
\SB{1.7}{2.9}{2.5}{sec::support}{Support Functions}
\draw (6.0,2.3) rectangle (8.6,3.1);
\SB{7.3}{2.9}{2.5}{sec::curvature}{Curvature}
\draw (1.2,1.4) rectangle (3.8,2.2);
\SB{2.5}{2.0}{1.6}{sec::complexity}{Complexity}
\draw (5.2,1.4) rectangle (7.8,2.2);
\SB{6.5}{2.0}{1.6}{sec::volume}{Volume}
\draw (3.1,0.5) rectangle (5.9,1.3);
\SB{4.5}{1.1}{0.7}{sec::zonotopes}{Zonotopes};

\draw (11.8,6.5) rectangle (15.2,7.3);
\node at (13.5,7.1) {\footnotesize \textbf{Section}~\ref{sec::reachable}};
\node at (13.5,6.7) {\footnotesize Reachable Sets};

\draw (10.2,5.2) rectangle (12.8,6.4);
\SBT{11.5}{6.2}{5.8}{sec::RPI}{Robust Positive}{5.4}{Invariant Set}
\draw (14.2,5.2) rectangle (16.8,6.4);
\SBT{15.5}{6.2}{5.8}{sec::RCI}{Robust Control}{5.4}{Invariant Set}

\draw (9.4,4.3) rectangle (12.0,5.1);
\SB{10.7}{4.9}{4.5}{sec::ControlSets}{Controllable Sets}
\draw (15.0,4.3) rectangle (17.6,5.1);
\SB{16.3}{4.9}{4.5}{sec::CLF}{CLFs};
\draw (9.4,3.4) rectangle (12.0,4.2);
\SB{10.7}{4.0}{3.6}{sec::CISets}{Control Invariance}
\draw (15.0,3.4) rectangle (17.6,4.2);
\SB{16.3}{4.0}{3.6}{sec::RCLF}{Robust CLFs};

\draw (10.2,2.1) rectangle (12.8,3.3);
\SBT{11.5}{3.1}{2.7}{sec::RCT}{Robust}{2.3}{Control Tubes}
\draw (14.2,2.1) rectangle (16.8,3.3);
\SBT{15.5}{3.1}{2.7}{sec::MPC}{Model}{2.3}{Predictive Control}

\draw (11.8,1.2) rectangle (15.2,2.0);
\node at (13.5,1.8) {\footnotesize \textbf{Section}~\ref{sec::RMPC}};
\node at (13.5,1.4) {\footnotesize Robust MPC};
\end{tikzpicture}
\end{center}
\caption{This article is structured into two key sections. Section~\ref{sec::PolyhedralComputing} provides a comprehensive survey of pivotal findings in polytopic computing. We delve into the intricate relationship between the geometric structure and computational complexity of parametric polyhedra, fostering a profound understanding of these concepts. Building on this insight, Section~\ref{sec::PolyhedralControl} explores modern polytopic control design methodologies. This section concludes with an in-depth discussion on modern polyhedral design approaches through convex optimization, encompassing techniques for computing invariant sets, control Lyapunov functions (CLFs), and robust model predictive control.}

\end{figure*}

\subsection*{Contributions}
The three main contributions of this survey are summarized below.
\begin{enumerate}

\item We provide a simple and concise overview of the overwhelming amount of literature on the geometry and algebra of polyhedra and polytopes. Our focus is on unifying aspects that are of principal importance to applied mathematicians and engineers who wish to {\bf  optimize sets using convex optimization} tools.

\item We show how to reformulate a long list of challenging polyhedral computing tasks for constrained linear control systems as convex optimization problems. This includes a review of recent methods suited for {\bf convex optimization-based control synthesis} via (robust) control invariant polytopes.

\item We present a comprehensive overview of {\bf convex optimization-based robust model predictive control} formulations, with a focus on the most significant advancements from the past two decades.

\end{enumerate}
Throughout, the focus is on polyhedral computing tools that hold paramount practical significance in applications. Our goal is to equip readers with a comprehensive toolbox of polyhedral methodologies tailored to address a broad spectrum of set-based control challenges.

\subsection*{Notation}
Throughout this article, the $n$-dimensional Euclidean space is denoted by $\mathbb R^n$. A set $X \subseteq \mathbb R^n$ is called convex if
\[
\forall x,x' \in X, \ \forall \vartheta \in [0,1], \quad \vartheta x + (1-\vartheta) x' \in X .
\]
Moreover, the interior of $X $ in $\mathbb R^n$ is denoted by $\interior{X}$. With $A \in \mathbb R^{m \times n}$ denoting a matrix, we define
\begin{eqnarray}
AX & \defeq & \{  Ax \ | x\in X \} \notag \\
\text{and} \quad X+X' & \defeq &  \{ x+x' \ | \ x\in X, \ x'\in X' \}. \notag
\end{eqnarray}
The latter operation is called the Minkowski sum of $X$ and $X'$, for arbitrary sets $X,X' \subseteq \mathbb R^n$. Finally, the convex hull of the column vectors of a given matrix $Y$ is denoted by
\[
\cvx{Y} \ \defeq \ \left\{ \ Y \lambda \ \middle| \ \sum_i \lambda_i = 1, \lambda \geq 0 \ \right\}.
\]
Here, inequalities involving vectors and matrices, including expressions such as ``$\lambda \geq 0$" in the aforementioned definition, are to be understood component-wise. A very similar convention is used for cones. Namely, if $Z$ is a given matrix, we use 
\[
\cone{Z} \ \defeq \ \left\{ \ Z \nu \ \middle| \ \nu \geq 0 \ \right\}
\]
to denote its associated cone. It can be interpreted as a homogenization of the convex hull of the columns of $Z$. We use $\one = (1,1,\ldots,1)^\tr \in \mathbb R^n$ to denote the vector in $\mathbb R^n$ whose components are all equal to $1$. The symbol $\mathbb 1$ denotes the $(n\times n)$ identity matrix. 

\section{Polyhedral Computing}
\label{sec::PolyhedralComputing}

In this section, we discuss the primary research on polyhedral computing. We focus here on practical classes of parametric polyhedra and scalable algorithms that are of fundamental importance in applications, such as set approximation, convex optimization, and control.

\subsection{Polyhedra}
\label{sec::polyhedra}
Let $F \in \mathbb R^{\nf \times n}$ be a matrix and $y \in \mathbb R^\nf$ a vector. The set of solutions of the linear inequality $F x \leq y$ is called a polyhedron. It is a convex set in $\mathbb R^n$, denoted by
\[
\Poly(F,y) \ = \ \{ x \in \mathbb R^n \mid F x \leq y \}.
\]
The set $\Poly(F,y)$ is non-empty if and only if there exists a vector
$z \in \mathbb R^\nf$ with
\[
z \geq 0, \ F^\tr z = 0, \quad \mathrm{and} \quad y^\tr z < 0.
\]
This result is known as Gale's theorem~\cite{Gale1960}, which in turn, can be interpreted as a consequence of Farkas' lemma~\cite{Farkas1902}. Another useful characterization of polyhedra is known under the name Minkowski-Weyl theorem~\cite{Minkowski1989,Weyl1934}, which states that $\Poly(F,y)$ is nonempty if and only if there exist matrices $Y$ and $Z$ such that
\[
\Poly(F,y) \ = \ \cvx{Y} + \cone{Z},
\]
recalling that $\cvx{Y}$ denotes the convex hull of the column vectors of $Y$ and $\cone{Z}$ the convex cone spanned by the column vectors of $Z$. This cone is called the recession cone of $\Poly(F,y)$ and turns out to satisfy
\[
\cone{Z} \ = \ \Poly(F,0)\;.
\]
The above relation can be found, for example, in~\cite[Prop.~1.12]{Ziegler1995} 
together with a more detailed discussion on recession cones and homogenization. 

Note that in the context of parameterizing convex sets, a principal motivation behind the study of polyhedra is their ability to approximate arbitrary compact convex sets. Namely, there exists for every $\epsilon > 0$ and any given compact convex set $X$ with 
non-empty interior, a polyhedron $\mathcal P(F,y)$ such that the Hausdorff distance between $F$ and $X$ is smaller than $\epsilon$. However, the smallest number $\nf$ of linear inequalities that are needed to represent such polyhedral approximations scales in general exponentially with the dimension $n$. For more accurate analysis and results regarding this number in dependence on $\epsilon$, $n$, and the curvature of $X$, we refer to~\cite{Schneider1981}.

\subsection{Representation}
\label{sec::representation}
The Minkowski-Weyl theorem, reviewed above, implies that polyhedra can either be represented in half-space format by storing $(F,y)$ or in vertex-cone format by storing $(Y,Z)$. Additionally, depending on the application, the following definitions and remarks are useful. 

\begin{itemize}

\item The halfspace representation $(F,y)$ of a polyhedron is called minimal if there is no other representation $\Poly(F,y) = \Poly(F',y')$, where $F'$ has fewer rows than $F$.

\item Similarly, the vertex-cone representation $(Y,Z)$ of a polyhedron is called minimal if there is no other representation $(Y',Z')$ with fewer columns.

\item The polyhedron $\Poly(F,y)$ contains $0$ in its interior if and only if $y > 0$. In this case, it is possible to re-scale the facet matrix, $F_i' \defeq F_i/y_i$, such that
\[
\Poly(F',\one) \ = \ \Poly(F,y)
\]
corresponds to a scaled representation.

\end{itemize}
Efficient polynomial run-time algorithms for finding minimal halfspace representations of polyhedra typically proceed by recursively detecting redundant half-space constraints. Therefore, these algorithms are also known under the name \textit{constraint reduction methods}; see~\cite{Tits2006} for a detailed discussion.

\subsection{Polytopes}
\label{sec::polytope}
A bounded polyhedron is called a polytope~\cite{Bronsted1983,Gruenbaum1967}. For a non-empty polyhedron $\Poly(F,y) \neq \varnothing$, the following statements are equivalent.
\begin{enumerate}

\item $\Poly(F,y)$ is a polytope.

\item There exists a matrix $Y$ such that
\[
\Poly(F,y) \ = \ \cvx{ Y }.
\]

\item We have $Fx \leq 0$ if and only if $x=0$.

\item There exists a non-negative matrix, $\Lambda \geq 0$, such that
\[
\Lambda F \ = \ [-\mathbb{1} \ \one].
\]

\item We have $0 \in \interior{\cvx{F^\tr}}$.

\end{enumerate}
Note that the equivalence of the first two statements is a direct consequence of the Minkowksi-Weyl theorem. Similarly, the third and fourth conditions are variants of Gale's theorem, which can also be proven by using Farkas' theorem of the alternative~\cite{Farkas1902}. Finally, the equivalence of the third and the last statement is an immediate consequence of Caratheordory's theorem~\cite{Caratheodory1911}.

\subsection{Linear Programming}
\label{sec::LP}
The analysis of convex polyhedra is closely related to the analysis of the generic linear programming problem, as introduced by Dantzig~\cite{Dantzig1963},
\begin{align}
\label{eq::LP}
\max_x \ c^\tr x \quad \mathrm{s.t.} \quad F x \leq y.
\end{align}
Namely, the feasible set of~\eqref{eq::LP} is the polyhedron $\Poly(F,y)$. Its objective value is called the support of $\Poly(F,y)$ in the direction $c \in \mathbb R^n$. The associated dual problem is given by
\begin{align}
\label{eq::DualLP}
\min_{z \geq 0} \ y^\tr z \quad \mathrm{s.t.} \quad F^\tr z = c
\end{align}
and the following weak and strong duality statements hold in full generality, without needing constraint qualifications~\cite{Rockafellar1970}.
\begin{enumerate}
\item For any $x \in \Poly(F,y)$ and any feasible point $z$ of~\eqref{eq::DualLP}, we have $c^\tr x \leq y^\tr z$.
\item If~\eqref{eq::LP} and~\eqref{eq::DualLP} are feasible, there exists a pair $(x^\star,z^\star)$ of feasible solutions with $c^\tr x^\star = y^\tr z^\star$.
\end{enumerate}

The primal and dual LP constitute convex optimization problems in the variables $x$ and $z$. However, the primal feasibility condition $Fx \leq y$ exhibits joint convexity in the variables $(x,y)$. This contrasts with the objective function of the dual LP, where a bilinear term arises. Specifically, the weak duality relation,
\begin{align}
\label{eq::wd}
c^\tr x \ \leq \ y^\tr z,
\end{align}
can be interpreted as a convex inequality in $(x,z)$ that fails to exhibit joint convexity in $(x,y,z)$. Indeed, as we shall elaborate further, the bilinear term $y^\tr z$ plays a pivotal role in posing challenges in convex optimization-based polyhedral computing. 

\subsection{Faces}
\label{sec::faces}
Faces can be used to analyze the structure of polyhedra and polytopes. Here, a face of a polyhedron $\Poly(F,y)$ is any set of the form
\[
\Poly(F,y) \cap \left\{ \ x \in \mathbb R^n \ \middle| \ a^\tr x = b \  \right\},
\]
where $(a,b) \in \mathbb R^n \times \mathbb R$ is such that we have $a^\tr x \leq b$ for all $x \in \Poly(F,y)$. The set $\varnothing$ is a trivial face, and all non-empty faces are called proper faces. The faces of dimension $0$, $1$, and $\mathrm{dim}(\Poly(F,y))-1$ are called vertices, edges, and facets of $\Poly(F,y)$. More general faces of dimension $k$ are called $k$-faces. Some basic but important consequences of this definition are summarized next.

\begin{itemize}
\item The proper faces of polyhedra are polyhedra.
\item The proper faces of polytopes are polytopes.
\item The intersections of faces of polyhedra are also faces.
\item The set of maximizers of~\eqref{eq::LP} is a face of $\Poly(F,y)$.
\item Every polyhedron has at most finitely many faces.
\item Every polytope is the convex hull of its vertices.
\end{itemize}
The first three statements are immediate consequences of the definition of faces. The fourth statement follows after setting $a = c$ and $b = c^\tr x^\star$ in the above definition with $x^\star$ denoting a maximizer of~\eqref{eq::LP}. The fifth statement follows from the fourth statement and the fact that strong duality holds for LPs; see~\cite[Thm.~2.14]{Fukuda2020}. A simple proof of the last statement based on Farkas' lemma can, for example, be found in~\cite[Prop.~2.2]{Ziegler1995}.

Faces are the building blocks used to describe the shape and structure of polytopes and polyhedra. To understand these shapes better, it is important to know the basic terms used to describe them.
\begin{enumerate}

\item The \textit{face lattice} of a polyhedron refers to the partially ordered set of its faces. In this context, the faces are ordered with respect to inclusion. In particular, the empty set is the minimal element of the lattice of every polyhedron and the polyhedron itself 
is its maximal element.x

\item The \textit{graph} of a polytope refers to its collection of vertices and edges, which form a simple undirected graph. And,

\item the \textit{$\mathsf{f}$-vector} of a polyhedron refers to the $(n+2)$-dimensional integer vector $(\mathsf{f}_{-1},\mathsf{f}_0, \ldots, \mathsf{f}_n)^\tr$. Here, one defines $\mathsf{f}_{-1} = 1$ and $\mathsf{f}_k$ denotes the number of the $k$-faces of the polyhedron for $0 \leq k \leq n$.
\end{enumerate}
Regarding results from the field of algebraic topology~\cite{Hatcher2002}, it is well-known that the $\mathsf{f}$-vector of a polytope satisfies the Euler-Poincar\'{e} equation
\[
\sum_{k=-1}^{n} (-1)^k \nf_k \ = \ 0.
\]
All other algebraic relations of the coefficients of the $\nf$-vector have been found and categorized by Dehn and Sommerville; see, for example,~\cite{Gruenbaum1967} or~\cite[Sect.~8.3]{Ziegler1995}. Moreover, a classical result on the combinatorial structure of polytopes, due to Balinski~\cite{Balinski1961}, states that the graph of $n$-dimensional polytopes is $n$-connected. In this context, however, it is important to be aware of the fact that analyzing the structural properties of general polytopes can be a difficult task. A prime example of this challenge is the enduring Hirsch conjecture, which posited an upper bound on the combinatorial diameter of polytopes,
i.e. the maximum shortest path between any two vertices of its graph, and remained unresolved for over half a century until its falsity was conclusively demonstrated in \cite{Santo2012}.  Nevertheless, as we will discuss below, graphs and face lattices of many practical classes of polytopes are much easier to analyze.

\subsection{Support Functions}
\label{sec::support}
The support function $\support{X}: \mathbb R^n \to \mathbb R \cup \{ \infty \}$ of an arbitrary closed set $X \subseteq \mathbb R^n$ is defined by 
\[
\forall c \in \mathbb R^n, \qquad \support{X}(c) \ \defeq \ \max_{x \in X} \, c^\tr x.
\]

It corresponds to the objective value obtained by solving~\eqref{eq::LP} for a given vector $c$. It is important to note, however, that due to the non-joint convexity of condition~\eqref{eq::wd} in $(x,y,z)$, the map
\[
(y,c) \quad \to \quad s[\Poly(F,y)](c)
\]
is in general not jointly convex in $(y,c)$, as discussed in~\cite[Sect.~3]{Rakovic2013}. This limitation underscores the restricted usefulness of support functions in analyzing polyhedra's properties as functions of their parameters. Nevertheless, the concept of support functions finds its strength in the following general principles, which are valid for all closed convex sets $X,X' \subseteq \mathbb R^n$ possessing non-empty interiors.
\begin{itemize}

\item The function $s[X]$ is positive homogeneous.

\item The function $s[X]$ is convex.

\item If $X \subseteq X'$ then $s[X] \leq s[X']$.

\item We have $s[X+X'] = s[X] + s[X']$.

\item If $s[X] = s[X']$, then $X = X'$.

\item $X$ is a polyhedron if and only if $s[X]$ is a piecewise affine function, defined on a finite polytopic partition of $\mathbb R^n$, 
i.e. a finite union of $n$-dimensional polyhedra with nonempty interior 
whose interiors are pairwise disjoint.  
\end{itemize}
The first four statements are valid for all closed sets $X$ and $X'$, as they directly follow from the definition of $s[X]$. For a detailed examination of the fifth statement, we refer to~\cite{Rockafellar1970}. Additionally,~\cite{Borelli2003} provides a proof for the sixth statement.

\subsection{Complexity}
\label{sec::complexity}
Because $\Poly(F,y)$ has only a finite number of faces, they can---at least in principle---be enumerated. For example, any index set $I \subseteq \{ 1, \ldots, \nf \}$ can be associated with a face of $\Poly(F,y)$, given by
\[
H_I(F,y) \ \defeq \ \left\{ \ x \in \Poly(F,y) \ \middle| \ \forall i \in I, \ F_i x = y_i \ \right\}.
\]
Since any proper face can be written in the form for at least one index set $I$ with cardinality $|I| \leq n$, it follows that $\Poly(F,y)$ has at most $\sum_{i=0}^n \binom{\nf}{i}$ proper faces. This bound is tight for simplices, which are defined as the convex hull of $n+1$ affinely independent vectors in $\mathbb R^n$. For other classes of polytopes, we might have $H_I(F,y) = \varnothing$ for many choices of $I$, but the total number of faces of many practical polyhedra and polytopes is, unfortunately, often very large. In this context, the three most important examples are:

\begin{enumerate}

\item any simplex in $\mathbb R^n$ has ``only" $n+1$ vertices and ``only" $n+1$ facets, but it still has
\[
\sum_{i=0}^n \binom{n+1}{i} \ = \ 2^{n+1}-1
\]
proper faces in total;

\item the unit cube $\{ x \in \mathbb R^n \mid \| x \|_\infty \leq 1 \}$ has ``only" $2n$ facets, but it has $2^n$ vertices and $3^n$ proper faces in total; and

\item the $1$-norm ball $\{ x \in \mathbb R^n \mid \| x \|_1 \leq 1 \}$ has ``only" $2n$ vertices, but it has $2^n$ facets and $3^n$ proper faces in total.
\end{enumerate}
Polytopes with a specified vertex count that exhibit the maximum number of facets can be constructed by determining the convex hull of distinct points along the moment curve, defined as
\[
t \to (t, t^2, \ldots, t^n)^\tr, \quad t \in \mathbb{R}.
\]
This approach results in cyclic polytopes, which are the key for proving McMullen's upper bound theorem~\cite{Mullen1971}; see also~\cite{Stanley1980}. This result implies that a polytope with $\nv$ vertices can have up to $\order{\nv^{\lfloor n/2 \rfloor}}$ facets. However, in practical applications, it is often desirable to consider polytopes that possess a balanced number of vertices and facets rather than those with an extreme facet or vertex count. For instance, convex polytopes that can be partitioned into a reasonable number of simplices represent a practical choice~\cite{Edmonds1970}. 
Other practical choices include some (but not all) classes of generalized permutohedra~\cite{Postnikov2006}, as well as various classes of random polytopes~\cite{Barany2007,Borgwardt1987,Spodarev2013}. For example, the number of expected facets of the convex hull of $v$ points, sampled uniformly from the unit sphere, scales as $\order{\log(\nv)^{n-1} \nv}$ for $\nv \to \infty$, see~\cite{Buchta1985} and~\cite[Thm.~1]{Newman2019}.

\subsection{Simplicity}
\label{sec::simplicity}
The regularity of solutions to parametric linear programs has been thoroughly investigated in the domains of convex optimization, as detailed in~\cite{Rockafellar1970}, and parametric optimization, as discussed in~\cite{Robinson1980}. Additionally, geometers are familiar with an equivalent concept of such regularity in the analysis of polytopes, which they refer to as simplicity, as outlined in~\cite{Gruenbaum1967}. In an attempt to unify these notions, we propose to call a proper face $H_I(F,y) \neq \varnothing$ of a polyhedron $\Poly(F,y)$ locally stable, if there exists an $\epsilon > 0$ such that
\[
H_I(F,y') \neq \varnothing,
\]
for all $y'$ with $\| y-y'\| \leq \epsilon$. Moreover, a polyhedron $\Poly(F,y)$ with minimal representation $(F,y)$ is called simple if all of its proper faces are locally stable. If $(F,y)$ is a minimal representation of a non-empty polytope $\Poly(F,y)$, then the following statements are equivalent.

\begin{itemize}

\item The polytope $\Poly(F,y)$ is simple.

\item All vertices of $\Poly(F,y)$ lie in exactly $n$ edges.

\item Every $k$-face of $\Poly(F,y)$ lies in $n-k$ facets, $k \geq 0$.

\item LP~\eqref{eq::DualLP} has for all $c \neq 0$ a unique minimizer $z^\star$.

\end{itemize}
The equivalence between the first and the last statement corresponds to a standard stability theorem for LPs that is well-known in the field of parametric optimization~\cite{Rockafellar1970}. The equivalence of the first two statements is discussed in~\cite{Gruenbaum1967}, see also~\cite{Shephard1963} or~\cite[Cor.~2]{Villanueva2024}. And, finally, the equivalence of the second and the third statement can be found in the literature on convex polytopes, see~\cite[Prop.~2.16]{Ziegler1995} or also~\cite[Prop.~4.10]{Fukuda2020}.

Numerous classes of simple polytopes hold significant practical importance. A pivotal reason for this lies in Perles' conjecture, which holds true for simple polytopes; it asserts that the face lattice of these polytopes is fully determined by their graph. This observation was initially established in \cite{Blind1987}, and subsequently, Kalai provided an exquisitely elegant proof \cite{Kalai1988}, which has since become the prevalent version in most geometry textbooks. From a computational standpoint, this finding is intriguing as a simple polytope with $\nv$ vertices possesses $\frac{n \nv}{2}$ edges, leading to a storage complexity for its graph that is at most of order $\order{n \nv}$. Consequently, merely storing the graph suffices to determine the entire face lattice. As we will discuss further below, this explains why numerous enumeration problems for certain practical classes of simple polytopes remain computationally tractable, even though the underlying algorithms may be NP-hard for general polytopes. Indeed, numerous counterexamples for non-simple polytopes, documented in \cite{Haase2002}, underscore the invalidity of Perles' conjecture for general polytopes.

\subsection{Polarity}
\label{sec::polarity}
Halfspace and vertex representation of polytopes are closely linked via the concept of polarity. Generally, if $X \subseteq \mathbb R^n$ is any set in $\mathbb R^n$, its polar set is defined as
\[
X^* \ = \ \left\{ \ z \in \mathbb R^n \ \middle| \ \sup_{x \in X} \ z^\tr x \leq 1 \ \right\}.
\]
This notion of polarity is relevant in the context of analyzing polyhedra. This is because the strong duality statement from Section~\ref{sec::LP} implies that the polar set of a polyhedron is itself a polyhedron. Moreover, if $F$ satisfies any of the equivalent conditions from Section~\ref{sec::polytope}, such that $\Poly(F,\one)$ is a polytope, then we have
\[
\Poly(F,\one)^* \ = \ \cvx{F^\tr}.
\]
This statement follows directly from the definition of polarity by using strong duality; see~\cite[Sect.~3.3]{Blanchini2008} for details. Consequently, the facet representation of a given polytope in $\mathbb R^n$, containing $0$ in its interior, can be used directly as a vertex representation of its associated polar set. By working with minimal representations, this result can be used to prove that the number of $k$-faces of $\Poly(F,\one)$, $0 \leq k \leq n-1$, equals the number of $(n-k-1)$-faces of its polar polytope~\cite[Sect.~2.3]{Ziegler1995}. And, in turn, by using the characterization of simple polytopes from the previous section, it follows that the polar set of a simple polytope of the form $\Poly(F,\one)$ is simplicial, meaning that all facets of $\cvx{F^\tr}$ are simplices.

Apart from the above explicit construction other statements about polar polyhedra and polytopes follow from more general results from convex analysis. For example, if $X$ is a closed convex set with $0 \in \interior{X}$, we have $(X^*)^* = X$, see~\cite{Rockafellar1970}. Obviously, this statement holds in particular for any polyhedron that contains $0$ in its interior. Moreover, for arbitrary sets $X_1,X_2 \subseteq \mathbb R^n$ the implication
\[
X_1 \subseteq X_2 \qquad \Longrightarrow \qquad X_1^* \supseteq X_2^*
\]
holds. This relation is of practical relevance, as it enables the transformation of polytopic outer approximations of a set into polytopic inner approximations of its polar set.

\subsection{Enumeration}
\label{sec::enumeration}
A matrix pair $(F,Z)$ is called a \textit{double-description pair} if it satisfies
\[
\Poly(F,0) \ = \ \cone{Z}.
\]
By the Minkowski-Weyl theorem for cones, there exists for every given $F$ a $Z$ such that this equation holds. In the algorithmic literature, the problem of computing a minimal $Z$ is known as the (vertex ray) enumeration problem~\cite{Fukuda1996}. The practical relevance of this problem stems from the fact that many problems in the field of computational geometry can be reduced to a vertex ray enumeration problem:

\begin{enumerate}

\item The \textit{conic facet enumeration problem} is about finding $F$ for a given $Z$ such $(F,Z)$ is a double description pair. Because Farkas' lemma implies that $(F,Z)$ is a double description pair if and only if $(Z^\tr,F^\tr)$ is a double description pair, this problem can be reduced the vertex ray enumeration problem~\cite[Page~92]{Fukuda1996}.

\item The \textit{convex hull problem} is about finding a facet representation of the convex hull, $\cvx{Y}$, of the columns that are stored in the given matrix $Y$. But this problem is equivalent to finding a facet representation of the cone
\[
\cone{ 
\begin{pmatrix}
Y \\
\one^\tr
\end{pmatrix}
 },
\]
also known as the homogenization of $\cvx{Y}$, see~\cite[Prop.~1.14]{Ziegler1995}. Namely, if a facet representation of the above cone is given by
\[
\cone{ 
\begin{pmatrix}
Y \\
\one^\tr
\end{pmatrix}
 }
\ = \ \Poly \left(
\left(
\begin{array}{cc}
F & -y \\
0 & -1
\end{array}
\right),\left(
\begin{array}{c}
0 \\
0
\end{array}
\right)
\right),
\]
then $(F,y)$ is the desired facet representation of $\cvx{Y}$. As such, the convex hull problem can be reduced to conic facet enumeration problem mentioned above and, hence, it can be solved by vertex ray enumeration.

\item The traditional \textit{vertex enumeration problem} for polytopes is polar to the convex hull problem and, hence, can be solved by vertex ray enumeration. More generally, we can find a minimal vertex-cone representation from the halfspace representation of an arbitrary polyhedron (and vice versa) by first constructing a double description of the recession cone and then applying the above strategy in combination with the Minkowski-Weyl theorem; see also~\cite{Avis1997}.

\item The \textit{Delaunay triangulation problem} is about computing a maximal triangulation of $\nv$ given points, $v_1,v_2, \ldots, v_\nv \in \mathbb R^n$ such that no point is (strictly) inside the $n$-dimensional circumsphere of each simplex. Because the regions of the Delaunay triangulation are the ``downwards-looking" faces of the lifted polytope,
\[
\cvx{ \left[
\left(
\begin{array}{c}
v_1 \\
\| v_1 \|^2
\end{array}
\right), \ \ldots, \ \left(
\begin{array}{c}
v_\nv \\
\| v_\nv \|^2
\end{array}
\right) \
\right] } \ \subseteq \ \mathbb R^{n+1},
\]
the Delaunay triangulation problem reduces to the convex hull problem~\cite{Brown1979}, see also~\cite[Page~693]{Graham1990} for an intuitive discussion. Hence, based on the above comments, this problem can be solved by vertex ray enumeration.

\item The \textit{post office problem} for $\nv$ given post-office sites, $x_1, \ldots, x_\nv \in \mathbb R^n$ is about dividing the world, $\mathbb R^n$, into $n$ postal regions $X_i$, such that the (Euclidean) distance from any point inside $X_i$ to its office site $x_i$ is smaller than to any other point $x_j$ with $j \neq i$. In the literature on computational geometry, the regions $S_i$ are also known under the name \textit{Voronoi regions}, see~\cite{Voronoi1908} and~\cite{Aurenhammer1991}. It is well-known that computing these regions is dual to the Delaunay triangulation problem~\cite{Brown1979}. Hence, the post office problem can be solved by vertex ray enumeration.
\end{enumerate}

Existing algorithms for solving the vertex ray enumeration problem can be divided into two classes, namely, incremental methods and pivoting methods. The former typically starts with an initial double description pair $(F_I,Y)$ for a subset $I$ of the row indices of $F$, which is then updated recursively by adding rows indices to the index set, $I \leftarrow I \cup \{ i \}$, while maintaining the double description. Examples of such methods are the original \textit{double description method}, which goes back to Motzkin~\cite{Motzkin1953}, and its polar version, the \textit{beneath-and-beyond} algorithm~\cite{Seidel1981}, see also~\cite{Fukuda1996} for a modern perspective and computational details.
In contrast, pivoting methods proceed by first locating a vertex (or extreme vertex ray) and recursively searching for neighboring vertices. This can be achieved by updating an adjacency graph~\cite{Avis1992}. As such, one could say that pivoting methods reverse the simplex algorithm in all possible ways, enumerating vertices and edges in a systematic way. An in-depth discussion of the differences between these methods and numerical comparisons can be found in~\cite{Avis1997}. Note that several libraries implement enumeration algorithms for polyhedra. The main ones are:
\begin{itemize}
\item {\tt qhull}~\cite{qhull}, which provides a C implementation of the divide-and-conquer Quickhull algorithm~\cite{Barber1996};
\item {\tt lrslib}~\cite{lrslib}, offering both single and multi-threaded C implementations of the reverse search vertex enumeration algorithm~\cite{Avis1992}  (see~\cite{Avis2000,Avis2018} for details); and
\item {\tt cddlib}~\cite{cddlib}, a C library that implements the Fourier-Motzkin double description method~\cite{Fukuda1996}.
\end{itemize}
Unfortunately, vertex enumeration algorithms and convex-hull computations, in general, exhibit an exponential worst-case complexity for arbitrary polyhedra and polytopes, as discussed in~\cite{Bremner1999} and~\cite{Klee1972}. However, for numerous practical classes of polytopes, these problems can be solved in polynomial time. For instance, in~\cite{Borgwardt1997,Borgwardt2007}, it is shown that the beneath-and-beyond algorithm for enumerating facets of random polytopes scales on average as $\order{\nv^2}$, where $\nv$ is the number of vertices. Similarly, double description methods for vertex enumeration of randomly generated polytopes with $\nf$ facets exhibit an expected runtime of $\order{\nf^2}$.
For simple $n$-dimensional polytopes with $\nf$ facets, the enumeration of their $\nv$ vertices using the reverse search method can be accomplished in $\order{\nv\nf n}$ time, with a storage complexity of $\order{\nv n}$~\cite{Avis2000}. Furthermore, simple polytopes can be efficiently computed from their graphs in polynomial time, as discussed in~\cite{Friedman2009}.
While the complete face lattice can theoretically be constructed from facet-vertex incidences, as outlined in~\cite{Kaibel2002}, such exhaustive enumerations are rarely required in practical applications. Moreover, for those interested in enumerative methods for computing the explicit solution map of parametric linear programs,~\cite{Borelli2003} and~\cite{Jones2005} provide valuable insights, with~\cite{Herceg2013} highlighting software implementations.

\subsection{Algebra}
\label{sec::algebra}
It is straightforward to demonstrate that the Minkowski sum of polytopes results in another polytope, and similarly, an analogous assertion holds true for polyhedra. Furthermore, it is well-established that each face of the Minkowski sum of two polyhedra corresponds to the Minkowski sum of faces from its constituent polyhedra~\cite[Thm.~15.1.1]{Gruenbaum1967}. However, unraveling the algebraic properties of these Minkowski sums proves to be a non-trivial task~\cite{Mullen1989}. From a computational standpoint, the situation is often even more challenging, as the problem of enumerating the facets of the Minkowski sum of polytopes belongs to the class of NP-hard problems, at least in its generality, as discussed in~\cite{Tiwary2008}. In general, the concept of homotheticity emerges as a valuable tool for analyzing Minkowski sums based on the following definitions.

\begin{itemize}

\item Two polytopes $P$ and $Q$ are called homothetic if there exists $a \in \mathbb R^n$ and $b \geq 0$ such that $P = \{ a \} + b Q$.

\item A polytope is called irreducible if it cannot be written as a Minkowski sum of two non-homothetic polytopes.

\end{itemize}
The following statement can be traced back to the work of Meyer, see, for example~\cite[Thm.~4]{Meyer1974}. It can be regarded as the foundation for analyzing the algebra of Minkowski sums:
\begin{itemize}
\item Every $n$-dimensional polytope with $\nf$ facets can be written as the Minkowski sum of at most $\nf-n$ irreducible polytopes.
\end{itemize}
Note that from a purely algebraic standpoint, the above statement indicates that Minkowski sums formally behave less like "sums" of a group (since no inverse elements exist) but more like "products" of a ring, which can be decomposed into irreducible factors---modulo homothetic equivalence. In fact, this observation has been formalized in a seminal article of McMullen on polytope algebras~\cite{Mullen1989}; see also~\cite{Mullen1993}.

\subsection{Configurations}
\label{sec::configurations}
The face configuration of a general parametric polyhedron $\Poly(F,y)$ refers to the collection of its non-empty faces, recalling that $\Poly(F,y)$ has only finitely many such faces, which can be enumerated by using index sets. Thus, let $\mathbb Y_I(F) \subseteq \mathbb R^\nf$ denote the set of parameters $y$ for which the face $H_I(F,y)$ is non-empty. For any given index set $I \subseteq \{ 1, \ldots, \nf \}$, this set can be written in the form
\begin{eqnarray}
\mathbb Y_I(F) &=& \left\{
\ y \in \mathbb R^\nf \ \middle|
\  \begin{array}{l}
\exists x \in \mathbb R^n: \ \forall i \in I, \\
F x \leq y, \ F_i x = y_i
\end{array} \
\right\}. \notag 
\end{eqnarray}
With $\mathrm{Im}(F) \defeq \{ Fx \mid x \in \mathbb R^n \}$ denoting the image set of $F$, statements about basic properties can be listed as follows.
\begin{itemize}
\item We have $\Poly(F,y) \neq \varnothing$ if and only if $y \in \mathbb Y_{\varnothing}(F)$.
\item We have $\mathrm{Im}(F) \subseteq \mathbb Y_I(F)$ for all $I \subseteq \{ 1, \ldots, \nf \}$.
\item We have $\mathbb Y_I(F) \subseteq \mathbb Y_{\varnothing}(F)$ for all $I \subseteq \{ 1, \ldots, \nf \}$.
\item The set $\mathbb Y_I(F)$ is for all $I \subseteq \{ 1, \ldots, \nf \}$ closed in $R^\nf$.
\item The set $\mathbb Y_I(F)$ is for all $I \subseteq \{ 1, \ldots, \nf \}$ a convex cone.
\item If $I \supseteq J$, then $\mathbb Y_{I}(F) \subseteq \mathbb Y_J(F)$.
\end{itemize}
The first four and the last statement are immediate consequences of the above definition, while a proof of the fifth statement can be found in~\cite[Lem.~1]{Villanueva2024}. Next, let
\[
\mathcal I(F,y) \ \defeq \ \{ \ I \subseteq \{ 1, \ldots, \nf \} \ \mid \ H_I(F,y) \neq \varnothing \ \}
\]
denote the collection of the index sets of the non-empty faces of the polyhedron $\Poly(F,y)$. It can be used to define the \textit{configuration cone},
\[
\mathbb C(F,y) \ \defeq \ \bigcap_{I \in \mathcal I(F,y)} \mathbb Y_I(F).
\]
In terms of being historically precise, it should be pointed out that the above configuration cone coincides with the Minkowski sum of the image set $\mathrm{Im}(F)$ and the closure of the so-called \textit{type cone}, which has originally been introduced in~\cite{Mullen1973}. As such, the conceptual ideas behind the introduction of configuration cones and type cones are the same, but working with configuration cones is usually more convenient in the context of practical computations, as we will discuss below. If $(F,y)$ is a minimal representation of a polytope $\Poly(F,y)$, then the following statements hold.

\begin{itemize}

\item If $y' \in \mathbb C(F,y)$ then $\mathbb C(F,y') \subseteq \mathbb C(F,y)$.

\item The set $\mathbb C(F,y)$ is a closed convex cone in $\mathbb R^\nf$.

\item $\Poly(F,y)$ is simple if and only if $\interior{\mathbb C(F,y)} \neq \varnothing$.

\item $\Poly(F,y)$ is irreducible $ \ \ \Longleftrightarrow \ \ $ $\mathrm{dim}(\mathbb C(F,y))=n+1$.

\end{itemize}
The first statement follows immediately from the definition of configuration cones. The second statement is a consequence of~\cite[Thm.~8]{Mullen1973}. The third statement follows, for example, from~\cite[Cor.~2]{Villanueva2024} and the fourth statement follows from the even more general observation that $\Poly(F,y)$ is a sum of $k$ irreducible polytopes if and only if $\mathrm{dim}(\mathbb C(F,y))=n+k$; see~\cite[Thm.~9]{Mullen1973}. As such, configuration cones (or type cones) can be regarded as the foundation of McMullen's polytope algebra~\cite{Mullen1989}. Additionally, it is worth noting that the third statement implies that the polytope $\Poly(F,y)$ is simple almost surely if $y$ is sampled randomly~\cite[Cor.~3]{Villanueva2024}; see also~\cite{Shephard1963} for related considerations.

Note that both type cones and configuration cones have been found to be of major practical interest, for example, in the contexts of learning polytopes from support function samples~\cite{Dosert2023} and robust control~\cite{Villanueva2024}. Apart from the above mentioned algebraic considerations, the reason why these cones are of interest is that there exist simple polynomial run-time algorithms for computing facet representations of the form
\begin{align}
\label{eq::CE}
\mathbb C(F,y) \ = \ \Poly(E,0),
\end{align}
with $E \in \mathbb R^{\nE \times \nf}$ denoting a facet matrix in the parameter space, compare \cite[Sect.~3.5]{Villanueva2024} and Section~\ref{sec::curvature} below. In order to design such an algorithm and to understand why the minimal representation dimension $\nE \in \mathbb N$ is a surprisingly small number for many practical classes of polyhedra and polytopes, however, one first needs to introduce the concept of vertex curvature that is discussed below.

\subsection{Curvature}
\label{sec::curvature}
The Gauss-Bonnet theorem, which relates the Gaussian curvature of a manifold to its Euler number is considered, by many differential geometers, one of the most beautiful and profound results in the theory of surfaces~\cite[Chap.~13]{Pressley2012}. Generalizations of the original results by Gauss and Bonnet can be found, for example, in~\cite{Allendoerfer1943} for polyhedra and in~\cite{Federer1959} for more general non-smooth manifolds. In its most basic variant, in the context of polytopes, this generalization states that the sum of the curvatures $\kappa_i(F,y)$ of the vertices of a polytope $\Poly(F,y)$ is equal to the surface area of the $(n-1)$-dimensional unit sphere, $\mathbb S_{n-1} \subseteq \mathbb R^n$; that is
\[
\sum_{i=1}^{\nv} \kappa_i(F,y) \ = \ |\mathbb S_{n-1}|,
\]
where $\nv$ denotes the number of  vertices of $\Poly(F,y)$. In this context, the curvature $\kappa_i(F,y)$ is defined as the surface area of the intersection of $\mathbb S_{n-1}$ with the normal cone of the $i$-th vertex of $\Poly(F,y)$. Equivalently, one could also say that the normal fan of a polytope---that is the collection of all normal cones of a polytope---is complete, which means that the union of all normal cones is equal to $\mathbb R^n$; see also~\cite[Sect.~7.1]{Ziegler1995}.


An important consequence of the above variant of the Gauss-Bonnet theorem is that representations of the configuration cone can be found by enumerating vertices. Namely, if $(F,y)$ is a minimal representation of a simple polytope $\Poly(F,y)$ and $\mathcal V(y)$ its collection of vertex indices, then its configuration cone coincides with its vertex configuration; that is,
\begin{align}
\mathbb C(F,y) \ = \ \bigcap_{I \in \mathcal V(y)} \ \mathbb Y_I(F).
\end{align}
This equation holds because the curvature of the vertices of $\Poly(F,y)$ coincides with the curvature of the vertices of $\Poly(F,y')$ for all $y' \in \interior{\mathbb C(F,y)}$; see also~\cite[Thm.~2]{Villanueva2024}. Let $F_I$ and $\mathbb{1}_I$ denote the matrices that consist of the rows of, respectively, $F$ and the unit matrix $\mathbb{1}$ that are in the index set $I$. Because the matrices
\[
V_I \ \defeq \ F_I^{-1} \mathbb{1}_I
\]
can be used to locate the vertices $v_I = V_I y$ of $\Poly(F,y)$ for all $I \in \mathcal V(y)$, we additionally find the explicit representation
\begin{eqnarray*}
\mathbb Y_I(F) &=& \Poly(F V_I - I,0) \\
\text{and} \qquad \mathbb C(F,y) &=& \bigcap_{I \in \mathcal V(y)} \Poly(F V_I-I,0),
\end{eqnarray*}
compare \cite[Sect~3.5]{Villanueva2024}. Hence, an upper bound on the number $\nE$ of rows of the minimal facet representation matrix $E$ in~\eqref{eq::CE} is given by
\[
\nE \ \leq \ \nv \cdot \nf
\]
and---still assuming that $\Poly(F,y)$ is a simple polytope with $\nf$ facets and $\nv$ vertices---this minimal representation can be found in polynomial run-time by using standard constraint reduction methods~\cite{Tits2006}. In fact, the above upper bound on the representation complexity of configuration cones can be further sharpened by enumerating the directed edges of $\Poly(F,y)$. For instance, if $\mathcal E(y)$ denotes the set of index pairs of the vertices that are linked by a directed edge, it is sufficient to enforce the homogeneous inequalities $F_j V_I y \leq y_j$ merely along the edge graph, with $\{ j \} = J \setminus \{ I \cap J \}$, for all $(I,J) \in \mathcal E(y)$. Consequently, $\nE$ is bounded above by the number of edges of $\Poly(F,y)$, recalling that the face lattice of simple polytopes is completely determined by their graph~\cite{Kalai1988}.

\subsection{Volume}
\label{sec::volume}
The volume of an $n$-dimensional compact set $X$ is defined, as usual, by the $n$-dimensional Lebesgue integral 
\[
\Vol(X) \ \defeq \ \int_X 1 \, \mathrm{d}x.
\]
For the special case that $\mathcal P(F, \overline y)$ is a simple polytope with minimal representation and $\mathbb C(F, \overline y)$ its configuration cone, it is possible to find a practical closed-form expression for the volume of the polytopes $\Poly(F,y)$ for $y \in \mathbb C(F,\overline y)$. This is due to the considerations in the previous section, which imply that we can compute linear vertex maps
\[
\Lambda_1, \ldots, \Lambda_\mathrm{s}: \mathbb R^{\nf} \to \mathbb R^{n \times n}
\]
such that for all $y \in \mathbb C(F,\overline y)$ with $y \geq 0$,
\[
\mathcal P(F,y) \ = \ \bigcup_{i=1}^\mathrm{s} \cvx{ [ 0, \Lambda_i y ] } 
\]
is a partition of $\mathcal P(F,y)$ into $\mathrm{s}$ simplices. An in-depth discussion of such simplicial complexes, including bounds on the number of simplices in dependence on the dimension and vertex numbers, can be found in~\cite{Edmonds1970}.  By adding the volumes of these simplices, we find that the $n$-th root of the volume of $\mathcal P(F,y)$ is given by the positive homogeneous expression
\[
\Lvol(y) \ \defeq \ \sqrt[n]{\Vol( \mathcal P(F,y) )} \ = \ \left( \frac{1}{n!} \sum_{i=1}^s \left| \mathrm{det}( \Lambda_i y ) \right| \right)^\frac{1}{n},
\]
where the second equation holds for all $y \in \mathbb C(F,\overline y)$. Here, the additional requirement, $y \geq 0$, is not needed anymore as the above expression is, by construction, invariant under offset shifts; that is,
\[
\Lvol(y) \ = \ \Lvol(y + Fx)
\]
for all $x \in \mathbb R^n$. Moreover, it follows from the Brunn-Minkowski inequality~\cite{Gardner2002} that $\Lvol$ is concave on $\mathbb C(F,\overline y)$; see also~\cite[Thm.~7.1.1]{Schneider2013}. For two non-homothetic polytopes $P(F,y)$ and $P(F,z)$ with $y,z \in \interior{\mathbb C(F,\overline y)}$, the strict version of the Brunn-Minkowski inequality even implies that we have
\[
\Lvol( \vartheta y + (1-\vartheta) z ) \ > \ \vartheta \Lvol(y) + (1-\vartheta) \Lvol(z)
\]
for all $\vartheta \in (0,1)$. Consequently, one could say that the function $\Lvol$ is---modulo homotheticity---strictly concave on the interior of the configuration cone $\mathbb C(F,\overline y)$. 

An overview of the huge amount of existing theoretical complexity results and algorithms for pre-computing the matrices $\Lambda_i$ via efficient triangulation can be found in the detailed textbook~\cite{DeLoera2010} as well as in the concise review by Lee and Santos~\cite[Ch.~16]{Toth2017}.

\subsection{Zonotopes}
\label{sec::zonotopes}
Zonotopes are a class of multi-symmetric polytopes that admit practical parameterizations. In the context of set computations, one such practical parameterization is given by
\[
\Zono(G,c,y) \ \defeq \ \left\{ \ G x + c \ \middle| \ -y \leq x \leq y \ \right\}.
\]
Here, $G \in \mathbb R^{n \times m}$ is a matrix whose $m$ columns are the so-called generating vectors of the zonotope, $c \in \mathbb R^n$ is the center of the polytope, and $y \in \mathbb R^{m}$ a non-negative parameter vector that can be used for scaling. Note that $\Zono(G,c,y)$ is non-empty if and only if $y \geq 0$. Some basic statements regarding the geometry of zonotopes are:
\begin{itemize}

\item Points and intervals are zonotopes.

\item Zonotopes are centrally-symmetric polytopes.

\item Every face of a zonotope is a zonotope.

\item A non-empty polytope is a zonotope if and only if all of its $2$-faces are point symmetric.

\end{itemize}
A proof for the last statement is provided in~\cite{Schneider1970}, see also~\cite[Theorem 3.5.2]{Schneider2013}, while the remaining three statements are straightforward consequences of this overarching result. However, it is noteworthy that when formulating assertions regarding symmetry, one must exercise caution. For instance, there exist polytopes in all dimensions $n \geq 4$, which, while not being zonotopes, have the property that all their facets are centrally symmetric~\cite{Mullen1970}.

One of the main reasons why zonotopes are of interest for practical computations is that they can be written as a 
Minkowski sum of a singleton and $m$ centered line segments
\[
\Zono(G,c,y) = \{c\} + [-y_1 G_1,y_1 G_1] + \ldots + [-y_m G_m,y_m G_m].
\]
An important consequence of this representation of zonotopes is that the Minkowski sum formula
\[
\Zono(G,c,y)+\Zono(G,d,z) \ = \ \Zono(G,c+d,y+z)
\]
holds for all $y,z \geq 0$ and all $c,d \in \mathbb R^n$. Moreover, yet another important consequence of such Minkowski representations is that the configuration cone of zonotopes does neither depend on the center $c$ nor on the scaling parameter $y \geq 0$. This is because the vertices of $\Zono(G,c,y)$ are located at the points
\[
v_i \ = \ G \cdot \mathrm{diag}(y) \cdot \sigma_i(G) + c,
\]
where $\sigma_i(G)$ denotes the column of the vertex sign matrix $\sigma(G) \in \{ -1,0,1 \}^{m \times \nv}$. 
It is then again a simple consequence of the above mentioned variant of the Gauss-Bonnet theorem that $\sigma(G)$ does neither depend on $c$ nor $y$ and, hence, is an invariant; see also~\cite[Chap.~7]{Ziegler1995}. In summary, this implies that a vertex representation of $\Zono(G,c,y)$ is given by
\[
\Zono(G,c,y) \ = \ c + \cvx{ \ G \cdot \mathrm{diag}(y) \cdot \sigma(G) \ }.
\]
for all $y \geq 0$ and all $c \in \mathbb R^n$. A worst-case optimal algorithm for computing $\sigma(G)$ can be found in~\cite{Edelsbrunner1986}. Additionally, reverse search algorithms for computing $\sigma(G)$ can be found in~\cite{Avis1996}, which perform well on many practical classes of zonotopes, such as, for example, zonotopes that are constructed by random line segment generation~\cite{Schneider2022}.

Next, regarding the computation of the $n$-th root of the $n$-dimensional volume of a zonotope, apart from using generic triangulations, one option is to use McMullen's formula~\cite{Shephard1974},
\begin{align*}
\ZLvol(y) & \defeq \ \sqrt[n]{\Vol( \mathcal Z(G,c,y) )} \\
& = \ 2 \left( \sum_{I \in \mathcal I} | \mathrm{det}( G_I \mathrm{diag}(y_I) ) | \right)^\frac{1}{n}
\end{align*}
for all $y \geq 0$, where $\mathcal I$ denotes the set of all ordered subsets of $\{ 1, \ldots, m \}$ that have cardinality $n$, and $G_I$ and $y_I$ the rows of $G$ and $y$ that are in the index set $I$. Due to the Minkowski-Brunn theorem~\cite{Gardner2002}, $\ZLvol$ is concave on $\mathbb R_+^m$, positive homogeneous, and strictly concave on the (relative) open simplex $\{ y > 0 \mid \one^\tr y = 1 \}$.

Note that one of the main reasons why zonotopes are of practical interest is that they are closed under Minkowski summation, as the above mentioned sum formula holds. Unfortunately, however, it is impossible to accurately approximate arbitrary centrally symmetric convex sets in $\mathbb R^n$ by using zonotopes for $n \geq 0$. For example, Schneider conclusively demonstrated in~\cite{Schneider1983} that for any dimension $n \geq 3$, one can find a certain approximation accuracy constant $\gamma_n > 1$. This constant has the property that no zonotope $\Zono$ can exist that satisfies
\[
\mathcal B_1^n \subset \Zono \subset \gamma_n B_1^n \qquad \text{for} \qquad \mathcal B_1^n \defeq \{ x \in \mathbb R^n \mid \|x\|_1 \leq 1 \}.
\]
Importantly, this finding holds true regardless of the number $m$ of line segments used to represent the zonotope $\Zono$. Specifically, for $n=3$, this constant is given by $\gamma_3 = \frac{3}{2}$, but the situation gets much worse in higher dimensional spaces, as we can find a divergent sequence of constant $\gamma_n$; that is, $\gamma_n \to \infty$ for $n \to \infty$; see also~\cite[Rem.~3.5.9]{Schneider2013}. 

Zonotopes have been the starting point for the construction of many other classes of convex sets and set parameterizations. For example, compact sets that can be written as an infinite Minkowski sum of line segments, convergent in the Hausdorff metric, are called zonoids, as surveyed in~\cite{Bolker1969}. The set of zonoids is a class of centrally symmetric sets that is closed under Minkowski addition and dilation~\cite[Sect.~3.5]{Schneider2013}, arising naturally in many applications~\cite{Goodey1993}. In contrast to more general centrally symmetric convex sets, zonoids can be approximated by zonotopes up to any given accuracy, as analyzed in~\cite{Bourgain1989}. Examples for other classes of convex and non-convex set parameterizations that \mbox{are---in} one or the other \mbox{way---based} on Minkowski sums are the class of \text{sums of zonotopes and a single ellipsoid}~\cite{Burton1976}, the class of \textit{constrained zonotopes}~\cite{Scott2016,Raghuraman2022}, the class of \textit{zonotope bundles}~\cite{Althoff2011}, the class of \textit{polynomial zonotopes}~\cite{Althoff2013}, as well as the class of \textit{constrained polynomial zonotopes}~\cite{Kochdumper2023}.

\section{Polyhedral Control Design}
\label{sec::PolyhedralControl}
Polyhedral computing methods are powerful tools for analyzing constrained linear systems, expressed as
\begin{align}
\label{eq::system}
x_{k+1} \ = \ A x_k + B u_k + w_k. \
\end{align}
Here, $x_k \in \mathbb X$, $u_k \in \mathbb U$, and $w_k \in \mathbb W$ model non-empty, compact, and convex constraints for the state, control, and disturbance variables, respectively. Given this linear recursion, together with its constraints and uncertainty model, the subsequent survey primarily focuses on polyhedral control design and set-based tools that solve a single convex optimization problem. However, we do occasionally discuss less scalable iterative procedures where they aid in understanding the limitations and trade-offs associated with different approaches. While primarily focused on scenarios where $A$ and $B$ are given matrices, the survey also encompasses methods that are applicable to multiplicative uncertainty. That is, the case where $(A,B)$ is an uncertain pair,
known to lie in some given compact and convex set.

\subsection{Reachable Sets}
\label{sec::reachable}
For the case that an affine feedback law~\cite{Kalman1960} is used, such that $u_k = u_\mathrm{r} + K (x_k - x_\mathrm{r})$, System~\eqref{eq::system} reduces to a constrained linear recursion with additive noise,
\begin{align}
\label{eq::linearRecursion}
\forall k \geq 0, \quad x_{k+1} \ = \ A x_k + w_k \quad \text{with} \quad w_k \in \mathbb W,
\end{align}
after redefining $A \leftarrow A + BK$ and $\mathbb W \leftarrow \mathbb W + \{ B u_\mathrm{r} - K x_\mathrm{r} \}$; see also~\cite{Hennet1995}. 
Under the assumption that the initial state $x_0 \in X_0$ is in a given set $X_0 \subseteq \mathbb R^{n_x}$, the set of all possible results for the state $x_k$ for arbitrary sequences $w_0,w_1,\ldots \in \mathbb W$ is denoted by $X_k \subseteq \mathbb R^{n_x}$. It is called the reachable set of~\eqref{eq::linearRecursion} at time $k$. Methods for computing or approximating $X_k$ can be classified into three categories, namely, \textit{exact methods}, \textit{indirect methods}, and \textit{recursive methods}.

\begin{itemize}

\item \textit{Exact methods} compute the set $X_k$ exactly or with high numerical precision. This is typically achieved by computing the Minkowski sums
\[
X_{k+1} = A X_k + \mathbb W.
\]
For example, if the initial set $X_0 = \Zono(G_0,c_0,\one)$ and the disturbance set $\mathbb W = \Zono(W,w,\one)$ are given zonotopes, the reachable sets $X_k = \Zono(G_k,c_k,\one)$ are found by computing~\cite{Girard2005,Girard2006}
\begin{align}
\forall k \in \mathbb N, \ \ G_{k+1} \ = \ [ A G_k, W ] \ \ \text{and} \ \ c_{k+1} \ = \ A c_k + w . \notag
\end{align}

Similarly, if $X_0 = \cvx{Y_0}$ and $\mathbb W = \cvx{W}$ are given polytopes with minimal vertex representation, any convex hull algorithm---see Section~\ref{sec::enumeration}---can be used to find minimal representations of the form 
\begin{align}
\cvx{Y_{k+1}} \ = \ \cvx{  (A Y_k) \otimes \one^\tr + \one^\tr \otimes W ] } \notag
\end{align}
recursively, for $k \in \mathbb N$, such that $Y_k$ is a minimal vertex representation of the polytope $X_k = \cvx{Y_k}$. Here, $\otimes$ denotes the Kronecker product, and the dimensions of the vectors $\one$ in the above expressions need to be adjusted such that all expressions match. 

\item \textit{Indirect methods} avoid the explicit computation of the sets $X_k$.  Instead, these methods aim at leveraging computationally tractable yet implicit conditions for representing inner and outer approximations as the solution of an optimization problem whose variables comprise worst-case disturbance scenarios. For example, if $F_1,F_2,\ldots,F_\nf \in \mathbb R^n$ are given facet normals, if $X_0 = \Poly(F,y_0)$ is a given polytope, and if $\mathbb W$ is polytope with given facet representation, too, then the LPs
\begin{eqnarray}
y_{k,j}^\star &=&  \max_{x_0,w,Y_{k,j}} F_j Y_{k,j} \notag \\
&  & \mathrm{s.t.} \quad
\left\{
\begin{array}{l}
Y_{k,j} = A^k x_0 + \sum_{i=0}^{k-1} A^i w_i \\
x_0 \in X_0, \ w_0,\ldots,w_{k-1} \in \mathbb W
\end{array}
\right. \notag 
\end{eqnarray}
for $j \in \{ 1, \ldots, \nf \}$ can be regarded as implicit representations of tight inner approximations of $X_k$, see~\cite{Guernic2010}. Namely, if $Y_{k,j}^\star$ and $y_{k,j}^\star$ denote maximizing state and objective value of the $j$th LP, polytopic inner and polyhedral outer approximations are given by
\[
\cvx{Y_k} \ \subseteq \ X_k \ \subseteq \ \Poly(F,y_k).
\]
Note that the outer approximation is tight in the sense that the facets of $\Poly(F,y_k)$ touch the sets $X_k$. 
Alternatively, one can solve an associated dual LP, if one wishes to express the support
\[
y_{k,j}^\star \ = \ \support{X_k}(F_j)
\]
as the optimal value of a minimization problem. Fully implicit representations, however, have the disadvantage that the number of implicit variables that are required for representing the tight inner and outer approximations of the sets $X_1,X_2,\ldots, X_N$ is of order $\order{N^2}$. This is in contrast to recursive methods that are discussed below.

\item \textit{Recursive methods} avoid the computation of disturbance scenarios. Instead, outer approximations of the reachable set $X_k$ are computed recursively. For example, if $\Poly(F,y_k) \supseteq X_k$ denotes a polytopic outer approximation of the reachable set at time $k$, an outer approximation of the form $\Poly(F,y_{k+1}) \supseteq X_{k+1}$ at the next time instance can be found by solving the LP
\begin{align}
\label{eq::polyrecursion}
\max_{x,y_{k+1}} \ \one^\tr y_{k+1} \quad \text{s.t.} \quad \left\{
\begin{array}{l}
\forall i \in \{ 1, \ldots, \nf \},\\
y_{k+1,i} \leq F_i A x_i + \overline w_i \\
F x_i \leq y_k,
\end{array}
\right.
\end{align}
where $\overline w_i = \support{\mathbb W}(F_i)$ denotes the support of the set $\mathbb W$ in the direction $F_i$. Variants of such recursive methods have also appeared~\cite{Girard2006}, where both inner and outer zonotopic approximations of the reachable set are computed. Moreover, recursive methods based on constrained zonotopes can be found~\cite{Scott2016}; see also~\cite{Raghuraman2021}.

\end{itemize}
Note that both indirect as well as recursive methods can be generalized for computing reachable sets of continuous time systems, as, for example, discussed in~\cite{Varaiya2000},~\cite{Villanueva2015}.

\subsection{Robust Positive Invariant Sets}
\label{sec::RPI}
In practical applications, one is often interested in computing the limit of the reachable sets $X_k$ of~\eqref{eq::linearRecursion} for $k \to \infty$,
\[
X_\infty \ \defeq \ \lim_{k \to \infty} X_k.
\]
Although one could define and analyze this (potentially unbounded) limit set for general linear systems, most of the existing methods for computing or approximating $X_\infty$ are based on the following assumption.
\begin{itemize}
\item[] \textbf{A1.} We assume that the eigenvalues of $A$ are all in the open unit disc and $\mathbb W$ is non-empty and compact.
\end{itemize}
If this assumption holds, $X_\infty$ is bounded, does not depend on the initial set $X_0$, and is given by the infinite Minkowski sum
\[
X_\infty \ = \ \sum_{k=0}^\infty A^k \mathbb W ,
\]
which converges absolutely in the Hausdorff metric. In terms of the structural properties of this set, if $A$ is invertible, $\mathbb W$ a compact convex set, and $\mathcal G$ the (multiplicative) cyclic group that is generated by $A$, it is known that $X_\infty$ belongs to the $\mathcal G$-invariant Minkowski class that is generated by $\mathbb W$~\cite[Sect.~3.5]{Schneider2013}. If $\mathbb W$ is a zonoid, it further follows that $X_\infty$ is a zonoid, too. From a computational perspective, however, this structural knowledge is certainly helpful but often also insufficient to compute the limit set. Therefore, in practice, one often uses positive invariant sets as surrogates for $X_\infty$. Here, a set $X \subseteq \mathbb R^{n}$ is called robust positive invariant (RPI) set for~\eqref{eq::linearRecursion} if
\[
AX + \mathbb W \subseteq X.
\]
The following statements about RPI sets are an immediate consequence of the above definition.
\begin{enumerate}
\item If $X$ is an RPI set, then $X^+ \defeq A X + \mathbb W$ is a RPI set.

\item If \textbf{A1} holds, then $X_\infty$ is an RPI set.

\item If \textbf{A1} holds and if $X$ is an RPI set, then $X_\infty \subseteq X$.

\item If \textbf{A1} holds, then $X_\infty$ is the minimal RPI set with 
respect to set inclusion.
\end{enumerate}
The first statement follows by multiplying the inclusion $AX + \mathbb W \subseteq X$ with $A$ from the left and adding $\mathbb W$ on both sides. Since this construction is such that $X^+ \subseteq X$, the first statement can be used to construct a monotonically decreasing sequence of RPI sets that converges to $X_\infty$ as long as \textbf{A1} holds; see also~\cite[Sect.~1]{Kolmanovsky1998}. The remaining three statements are thus a simple consequence of the first statement; see also~\cite[Chap.~4]{Blanchini2008} for details.

Historically, the first necessary and sufficient characterization of positive invariant polytopes can be traced back to the work of Bitsoris~\cite{Bitsoris1988}, from 1988, who proved that the polytope $\Poly(F,y)$ is positive invariant if and only if there exists a componentwise non-negative matrix $Z$ such that
\[
Z F = F A, \ Z \geq 0, \quad \text{and} \quad  Zy = y
\]
   Due to the presence of the bilinear term $Z y$, however, this condition is, unfortunately, not jointly convex in $(y,Z)$. Thus, 
   while this condition is well-suited for invariance verification, it
   is not for practical convex optimization-based computation of invariant polytopes.. Bitsoris also proposed an equivalent polar condition for vertex representation, which is, however, non-convex, too. Other older articles on analyzing and computing invariant polyhedra are can be found in~\cite{Genesio1985},~\cite{Vassilaki1989},~\cite{Bitsoris1999},~\cite{Stoican2015}  and~\cite{Dorea1999}.

A different approach that computes an RPI set of a priori specified, fixed complexity was developed by Trodden ~\cite{Trodden2016}. Namely, assuming that \textbf{A1} holds, Trodden showed that there exists an RPI set of the form $\Poly(F,y)$, with given facets matrix $F$ and $\overline w_i = \support{\mathbb W}(F_i)$, if and only if the LP
\begin{eqnarray}
\label{eq::Trodden}
\max_{x,y,y^+} \ \one^\tr y^+ \quad \text{s.t.} \quad \left\{
\begin{array}{l}
\forall i \in \{ 1, \ldots, \nf \},\\
y_i^+ \leq F_i A x_i + \overline w_i \\
F x_i \leq y, \ y^+ = y
\end{array}
\right.
\end{eqnarray}
admits a maximizer. And, in this case, denoting the maximizer for $y = y^+$ by $y^\star$, $\Poly(F,y^\star)$ is the smallest approximation, for the given $F$, to the minimal robust positive invariant polytope
in the inclusion sense. Moreover, reversely, if there exists no such polytopic RPI set, the above LP is unbounded; see~\cite[Thm.~4]{Trodden2016}. Regarding the basic intuition behind this construction, first note the above LP is obtained from~\eqref{eq::polyrecursion} by adding the fix-point constraint $y = y^+$. The fact one can jointly maximize over $y$ and $y^+$ in order to find a minimal RPI polytope is, however, not obvious at all. Here, one needs to use the monotonicity and positive homogeneity properties of the support function of the set $A \Poly(F,y)$ as a (non-convex) function of $y$; see~\cite{Rakovic2013}. A more recent discussion of the corresponding non-convex analysis techniques has also recently appeared in~\cite{Rakovic2024a}. Moreover, other related methods for computing RPI sets can be found in~\cite{Rubin2018,Rubin2018a}.

Finally, regarding finite time-horizon based approximations of RPI sets, one rather basic strategy proceeds by first approximating the set $X_\infty$ by the finite sum
\[
X_N \ = \ \sum_{k=0}^{N-1} A^k \mathbb W \ \approx \ X_\infty,
\]
recalling that $X_N$ is the $N$-step reachable set. Of course, in order to use this approximation, one needs to introduce a remainder bound. However, if $N$ is sufficiently large and \textbf{A1} holds, one can find an $\alpha \in [0,1)$ such that
\[
A^N \mathbb W \ \subseteq \ \alpha \mathbb W.
\]
Using the properties of geometric sums, it is not difficult to check that this condition is sufficient to guarantee that
\[
X = \frac{1}{1-\alpha} X_N
\]
is an RPI set. The details of this argument can, for example, be found in~\cite{Rakovic2005}, where also related numerical procedures for computing $\alpha$ are proposed. In this context, the set $X_N$ can either be computed exactly as in~\cite{Rakovic2005} or represented implicitly by using the indirect approach that has been detailed in the previous section. In fact, polytopic computing methods based on such implicit approach can be found~\cite{Rakovic2023}.

If $\mathbb W$ is a polytope, $X_N$ is a polytope and, hence, $(1-\alpha)^{-1} X_N$ is also a polytope. As such, the above strategy certainly yields invariant polytopic outer approximations of $X_\infty$. In~\cite[Sect.~IV]{Rakovic2005}, an algorithm is provided which, for every $\epsilon > 0$, one can compute $\alpha$ and $N$ such that \mbox{$(1-\alpha)^{-1} X_N$} is an $\epsilon$-outer-approximation of $X_\infty$. This means that arbitrarily accurate polytopic approximations of the minimal RPI set can be computed using this method. However, for small $\epsilon$, $N$ is typically large, resulting in a high computational cost to compute the corresponding sets $X_N$.

\subsection{Controllable Sets}
\label{sec::ControlSets}
The $N$-step controllable set $S_N$ is defined as the set of initial states $x_0$ for which there exists a control sequence $u_0,u_1,\ldots, u_{N-1} \in \mathbb U$ such that we have $x_N = 0$, where
\[
x_{k+1} = A x_k + B u_k, \quad k \in \{ 0, 1, \ldots, N-1 \}.
\]
 In this context, we typically assume $0 \in \mathbb U$, without loss of generality. Clearly, with $A^\dagger$ denoting the pseudo-inverse of $A$, computing $S_N$ is equivalent to computing $X_N$ after reversing time by re-defining $A \leftarrow A^\dagger$ and $\mathbb W \leftarrow A^\dagger (-B) \mathbb U$. In other words, we have
\[
S_N \ = \ \ker(A) + \sum_{i=0}^{N-1} \left( A^\dagger \right)^{k+1} (-B) \mathbb U,
\]
see also~\cite[Sect.~6.1.1]{Blanchini2008}. Consequently, if no state constraints are present, one could argue that computing $S_N$ is analogous to the computation of reachable sets, and, at least in principle, the same computational procedures can be used. The situation, however, becomes much more involved if state constraints, $\mathbb X \subseteq \mathbb R^n$, are present. In this case, starting from $S_0 = \{ 0 \}$, $S_N$ needs to be computed via the backward recursion
\begin{equation}\label{eq::kcontrollable}
S_{k+1} \ = \ \left( \ker(A) + A^\dagger S_k + A^\dagger (-B) \mathbb U \right) \cap \mathbb X,
\end{equation}
which is much harder to implement, as not only Minkowski sums but also intersections need to be computed. Of course, if $\mathbb U$ and $\mathbb X$ are polyhedra, the sets $S_k$ are polyhedra, too. However, even if no state constraints are present, the corresponding exact computations tend to become computationally demanding already for very small $k$, recalling that computing Minkowski sums of polyhedra is NP-hard in general~\cite{Tiwary2008}.

Because computing controllable sets exactly is difficult, most methods focus on computing inner approximations of $S_N$. In this context, one of the most profound results has been established in the 1980s by Gutman and Cwickel~\cite{Gutman1986,Gutman1987}. Namely, under the assumption that $\mathbb U$ and $\mathbb X$ are convex, their celebrated vertex control theorem states that two polytopes $S = \cvx{Y}$ and $S^+ = \cvx{Y^+}$, both given in vertex representation, satisfy
\begin{align}
\label{eq::CT}
S^+ \ \subseteq \ \left( \ker(A) + A^\dagger S + A^\dagger (-B) \mathbb U  \right) \cap \mathbb X
\end{align}
if and only if there exist vertex controls $U_i \in \mathbb U$ and a componentwise non-negative matrix $Z$ such that
\[
A Y^+ + B U = Y Z, \ Z \geq 0, \ Z^\tr \one \leq  \one, \ U_i \in \mathbb U, \ Y_i^+ \in \mathbb X;
\]
see also~\cite[Thm.~4.44]{Blanchini2008}. Unfortunately, similar to Bitsoris' condition for RPI sets, the above condition is not jointly convex in $(U,Y,Y^+,Z)$ due to the occurrence of the bilinear term, $YZ$. 
However, the above conditions can be turned into convex conditions if one restricts oneself to suitable classes of polytopes. For example, the following conditions are entirely convex.

\begin{enumerate}

\item Let $F$ be a facet matrix of a simple polytope with vertex matrices $V_1,\ldots,V_\nv$ and let $\Poly(E,0)$ be the facet representation of its associated configuration cone, such that
\[
\Poly(F,y) \ = \ \cvx{ V_1 y, \ldots, V_{\nv} y }
\]
for all $y$ with $E y \leq 0$, see Section~\ref{sec::curvature}. Then the configuration-constrained polytopes \mbox{$S = \Poly(F,y)$} and $S^+ = \Poly(F,y^+)$, with $E y \leq 0$ and $E y^+ \leq 0$, satisfy~\eqref{eq::CT} if and only if there exist vertex control inputs $U_i \in \mathbb U$ with
\begin{align}
\label{eq::CCCT}
\begin{array}{l}
F V_i y + F B U_i \leq y^+, \\
U_i \in \mathbb U, \ V_i y^+ \in \mathbb X, \ E y \leq 0, \ E y^+ \leq 0.
\end{array}
\end{align}
for all $i \in \{ 1, \ldots, \nv \}$. Note that this condition is jointly convex in $(U,y,y^+)$), whenever $\mathbb X$ and $\mathbb U$ are convex.

\item Let $G$ be a generator matrix and $\sigma(G)$ its associated vertex sign matrix. The zonotopes $S = \mathcal \Zono(G,c,y)$ and $S^+ = \mathcal \Zono(G,c^+,y^+)$ satisfy~\eqref{eq::CT} if and only if there exist vertex controls $U_i$ and auxiliary variables $x_i^+$ with
\begin{align}
\label{eq::ZCT}
\begin{array}{l}
G x_i^+ + c^+ = B U_i \\
\phantom{G x_i^+ + c^+ = } + A \left( G \cdot \mathrm{diag}(y) \cdot \sigma_i(G) + c \right) \\
-y^+ \leq x_i^+ \leq y^+, \ U_i \in \mathbb U, \
y \geq 0, \ y^+ \geq 0, \\
G \cdot \mathrm{diag}(y^+) \cdot \sigma_i(G) + c^+ \in \mathbb X
\end{array}
\end{align}
for all $i \in \{ 1, \ldots, \nv \}$. This condition, too, is convex in $(U,y,y^+,x^+)$ whenever $\mathbb X$ and $\mathbb U$ are convex.

\end{enumerate}
The above statements can be viewed as immediate consequences of Gutman's and Cwickel's vertex control theorem; see also~\cite[Cor.~4]{Villanueva2024} for details. Note that all of the above conditions refer to the computation of vertex control inputs. Such inputs however, can be interpolated to construct a nonlinear feedback control law; see, for example,~\cite{Nguyen2013} for details. 

\subsection{Control Invariant Sets}
\label{sec::CISets}
A set $S$ is called control invariant, if it satisfies~\eqref{eq::CT} for $S^+ = S$ and $\lambda$-control invariant, $\lambda \in [0,1]$ if it satisfies~\eqref{eq::CT} for $S^+ = \lambda S$. The following statements are immediate consequences of this definition, recalling that we assume $0 \in \mathbb U$ in this context.
\begin{enumerate}

\item The sets $S_k$ in~\eqref{eq::kcontrollable} are control invariant for all $k \in \mathbb N$.

\item If $\mathbb X$ is closed and bounded, the  limit of the controllable sets,
\[
S_\infty = \lim_{k \to \infty} S_k,
\]
exists. Moreover, every control invariant set $S$ satisfies $S \subseteq S_\infty$. Hence, $S_\infty$ is the maximal control invariant set with respect to set inclusion.
\end{enumerate}
The first statement holds as we have $S_{k+1} \supseteq S_k$ by construction. And, if $\mathbb X$ is closed bounded, the sequence $S_k$, $k \in \mathbb N$, is monotonically increasing and bounded, which implies the second statement; see also~\cite{Borelli2017} for details.

The relevance of the first statement is that any $k$-step controllable set $S_k$ can be used as a control invariant set. Alternatively, one could also solve convex optimization problems, for example, 
\begin{eqnarray}
\max_{y,y^+,U} \ \Lvol(y) \quad &\text{s.t.}& \quad
\left\{
\begin{array}{l}
 \eqref{eq::CCCT} \ \text{and} \\
\ y^+ = \lambda y \notag \\
\end{array}
\right.
\end{eqnarray}
or, alternatively,
\begin{eqnarray}
\max_{c,c^+,y,y^+,z^+,U} \ \ZLvol(y) \quad &\text{s.t.}& \quad \left\{
\begin{array}{l}
\eqref{eq::ZCT} \ \text{and} \\
y^+ = \lambda y, \\
c^+ = c = 0
\end{array}
\right. \notag
\end{eqnarray}
to compute $\lambda$-control invariant polytopes $\mathcal P(F,y)$ or zonotopes $\mathcal Z(G,c,y)$ with maximal volume, which approximate $S_\infty$ from inside. Other methods for computing maximal control invariant sets proceed by using iterative polytopic approximation methods, as worked out in much detail in~\cite{Athanasopoulos2013}; see also~\cite{Anevlavis2019}. Similar methods can also be found in the textbooks~\cite{Blanchini2008} and~\cite{Borelli2017}, which also provide a more complete overview of such iterative strategies for computing polytopic and non-polytopic control invariant sets.

\subsection{Robust Control Tubes}
\label{sec::RCT}
A sequence of sets, $X = (X_0,X_1, \ldots, X_N)$, is called a robust control tube (RCT) if
\[
\forall x \in X_k, \ \exists u \in \mathbb U: \ \forall w \in \mathbb W, \quad A x + Bu + w \in X_{k+1}
\]
for all $k \in \{ 0,1, \ldots, N-1 \}$. Here, we also formally allow the case $N = \infty$ for tubes that evolve on an infinite time horizon. Moreover, if state constraints $\mathbb X$ are present, we say that $X$ is a feasible RCT if $X_k \subseteq \mathbb X$ for elements $X_k$ of the sequence $X$. In the literature, RCTs are sometimes also called robust forward invariant tubes. RCTs are uniquely characterized by Kerrigan's condition, which he published for the first time in his PhD thesis~\cite[Section 2.10.2]{Kerrigan2000}. It states that a sequence of non-empty sets $X$ is an RCT if and only if we can find a sequence of non-empty intermediate sets $S_0,S_1,\ldots, S_{N-1}$, $S_k \neq \varnothing$, such that\footnote{In order to be historically precise, it should be pointed out that Kerrigan's original condition in~\cite{Kerrigan2000} used the notion of Pontryagin differences rather than intermediate sets. Apart from this, it should be mentioned that similar but more specialized conditions for RCTs have been around in the control literature of the 20th century. For example, A.B.~Kurzhanski and Valyi use a similar condition in the context of their ellipsoidal calculus~\cite{Kurzhanski1997}.}
\begin{eqnarray}
\label{eq::Kerrigan}
A X_k \subseteq S_k + (-B)\mathbb{U}  \qquad \text{and} \qquad S_k + \mathbb{W}  \subseteq X_{k+1}
\end{eqnarray}
for all $k \in \{ 0,1, \ldots, N-1 \}$.  The implications of this general condition are extensive and far-reaching. For example, if $\mathbb W$ and $\mathbb U$ are closed and convex, the subsequent assertions hold true.

\begin{enumerate}

\item The element-wise convex hull of an RCT is an RCT.

\item If $X$ and $X'$ are element-wise convex RCTs, then
\[
\vartheta X + (1-\vartheta)X' = \left( \ \ldots, \ \vartheta X_k + (1-\vartheta) X_k', \  \ldots \ \right)
\]
is an RCT for all $\vartheta \in [0,1]$.
\end{enumerate}
The first statement is a straightforward consequence of applying the convex hull operation to both sides of~\eqref{eq::Kerrigan}. The significance of this statement lies in the fact that when $\mathbb U$, $\mathbb W$, and $\mathbb X$ are all closed and convex, we can limit our investigation to element-wise convex RCTs. This is because the element-wise convex hull of a feasible RCT remains a feasible RCT. Furthermore, if the intermediate sets  $S_k$ and $S_k'$ satisfy~\eqref{eq::Kerrigan} for element-wise convex RCTs $X$ and $X'$, respectively, then the second statement can be verified by showing that the intermediate sets $\vartheta S_k + (1-\vartheta)S_k'$ satisfy the same condition for the convex combination $\vartheta X + (1-\vartheta) X'$. Hence, Kerrigan's condition establishes a pivotal connection between control theory and Brunn-Minkowski theory.

However, it is unfortunate to note that this algebraic convexity structure inherent in RCTs does not directly translate into practical convex conditions for general polyhedra of the form $\Poly(F,y)$. 
Consequently, current methods for analyzing and optimizing polyhedral or polytopic RCTs using convex analysis techniques often impose limitations, either by restricting the class of polytopes considered or by constraining the class of feedback laws that can generate the RCT, or both. Specifically, given closed and convex sets X and U, the existing convex parameterizations of RCTs can be categorized as follows.

\begin{itemize}

\item \textit{Scenario RCT} parameterizations~\cite{Kerrigan2004} enumerate all extreme uncertainty scenarios assuming that the sets
\[
X_0 = \cvx{Y} \quad \text{and} \quad \mathbb W = \cvx{W}
\]
are polytopes with a given vertex representation. The elements of the set of extreme scenarios,
\[
\mathcal S_N \defeq \{ Y_1, \ldots, Y_{\nv_x} \} \times \{ W_1, \ldots, W_{\nv_w} \}^N,
\]
are sequences of the form
\[
(x_0^i,w_0^i,\ldots,w_{N-1}^i) \in \mathcal S_N.
\]
They are enumerated by the index $i \in \mathcal I_N$, with $\mathcal I_N$ denoting the associated index set. The sequence
\[
X_k \ = \ \cvx{ \ \{ x_k^i \mid i \in \mathcal I_N \} \ }
\]
for $k \in \{0,\ldots N-1\}$
is then a feasible RCT if and only if there exists control inputs $u_{s(k,i)}$ such that
\begin{align}
\label{eq::ScenarioRCT}
\begin{array}{l}
x_{k+1}^i = A x_k^i + B u_{s(k,i)} + w_k^i, \\
x_k^i \in \mathbb X \quad \text{and} \quad u_{s(k,i)} \in \mathbb U
\end{array}
\end{align}
for all $i \in \mathcal I_N$. In this context, we enumerate the control vectors by the index $s(k,i) \in \mathbb N$ that satisfies $s(k,i) = s(k,j)$ whenever
\[
(x_0^i,w_0^i,\ldots,w_{k-1}^i) \ = \ (x_0^j,w_0^j,\ldots,w_{k-1}^j)
\]
and $s(k,i) \neq s(k,j)$ otherwise, for all $i,j \in \mathcal I_N$ with $i \neq j$; see~\cite{Kerrigan2004} and~\cite{Diehl2007}. Recalling our assumption that $\mathbb X$ and $\mathbb U$ are closed convex sets, this characterization of polytopic RCTs is convex. One significant drawback of this approach, however, is that the size of the index set $\mathcal I_N$ increases exponentially with the prediction horizon, $|\mathcal I_N| = \nv_x \nv_w^N$, which poses a significant computational challenge~\cite{Hubner2020}.

\item \textit{Separable Feedback RCT} parameterizations attempt to reduce the complexity of the above scenario enumeration approach by introducing a separable feedback structure. This means that any initial state $x_0$ and disturbance sequence $(w_0,\ldots,w_{k-1})$ is mapped to a control input
\[
u_k \ = \ \mu_{k,0}(x_0) + \sum_{j=0}^{k-1} \mu_{k,j+1}(w_j)
\]
at time $k$. Due to the vertex control theorem by Gutman and Cwickel~\cite{Gutman1986,Gutman1987} it is sufficient to parameterize the separable vertex inputs
\begin{align*}
& U_{k,0} \ \defeq \ \left[ \ \mu_{k,0}(Y_{0,0,1}), \ \ldots, \ \mu_{k,0}(Y_{0,0,\nv_x}) \ \right], \\
& U_{k,j} \ \defeq \ \left[ \ \mu_{k,j}(W_1), \ \ldots, \ \mu_{k,j}(W_{\nv_w}) \ \right].
\end{align*}
Here, $Y_{0,0,i} \in \mathbb R^{n_x}$ denotes the $i$-th vertex of the initial polytope $X_0 = \cvx{Y_{0,0}}$ and $W_i \in \mathbb R^{n_x}$ denotes the $i$-th vertex of the uncertainty set $\mathbb W = \cvx{W}$. Next, due to the above separable superposition structure of the feedback law, the parametric polytopes
\[
X_k \ = \ \sum_{j=0}^k \cvx{Y_{k,j}}
\]
are a feasible RCT, for $\mathbb X = \Poly(F,\overline y)$, $\mathbb U = \Poly(G,\overline u)$. and $\mathbb W = \cvx{W}$, if there exist separable vertex input matrices $U_{k,j}$ and vectors $y_{k,j},v_{k,j}$ with
\begin{align}
\label{eq::SeparableRCT}
\begin{array}{l}
Y_{k+1,0} = A Y_{k,0} + B U_{k,0}, \\
Y_{k+1,j} = A Y_{k,j} + B U_{k,j}, \ \ Y_{j,j} = W, \\
F Y_{k,0} \leq y_{k,0} \one^\tr, \ G U_{k,0} \leq v_{k,0}\one^\tr, \\
F Y_{k,j} \leq y_{k,j}\one^\tr, \ \ G U_{k,j} \leq v_{k,j}\one^\tr, \\
\overset{k}{\underset{j=0}{\sum}} y_{k,j} \leq \overline y, \ \ \overset{k}{\underset{j=0}{\sum}} v_{k,j} \leq \overline u.
\end{array}
\end{align}
for all $k \in \{ 0, \ldots, N \}$ and all $j \in \{ 1, \ldots, k \}$. Note that this parameterization has complexity $\order{N^2}$ with respect to the time horizon, and the above condition is convex; see also~\cite{Rakovic2011FP,Rakovic2012a} for implementation details. One should be clear in mind, however, that the above condition is merely sufficient for $X$ to be a feasible RCT but, in general, not necessary. This is because non-separable feedback laws may exist that generate tighter RCTs. Earlier but less general variants of the above method for affine disturbance feedback parameterizations---which are also separable---can be found in~\cite{Goulart2006}.

\item \textit{Rigid RCT} parameterizations~\cite{Langson2004,Mayne2005} assume that an asymptotically stabilizing linear control gain $K$ and an associated RPI set $\mathrm X_\mathrm{s} = \Poly(F,y_\mathrm{s})$ for the pre-stablized system,
\[
(A+B K) \Poly( F, y_\mathrm{s} ) + \mathbb W  \subseteq \mathcal \Poly(F,y_\mathrm{s}),
\]
are given, recalling that, for any such given $K$, an optimal $y_\mathrm{s}$ can be found by solving Trodden's LP~\eqref{eq::Trodden} with $A \leftarrow A + BK$. Next, the sequence of polytopes
\[
X_k = \{ x_k \} + \Poly(F,y_\mathrm{s})
\]
is---by construction---a feasible RCT as long as there exist central inputs $u_k$ such that
\begin{align}
\label{eq::RigidRCT}
\begin{array}{l}
x_{k+1} = A x_k + B u_k, \\
\{ x_k \} + X_\mathrm{s} \subseteq \mathbb X \quad \text{and} \quad \{ u_k \} + K X_\mathrm{s} \subseteq \mathbb U,
\end{array}
\end{align}
for $k \in \mathbb N$; see, for example,~\cite{Mayne2005} for details. Note that this approach parameterizes both the sequence of feedback laws, $\mu_k(x) = u_k + K(x-x_k)$ with $K$ being pre-computed, as well as the sequences of polytopes, $X_k = \{ x_k \} + \mathcal P(F,y_\mathrm{s})$ with $y_\mathrm{s}$ being pre-computed.

\item \textit{Homothetic RCT} parameterizations~\cite{Rakovic2012,Rakovic2013} are also based on pre-computing $K$ as well as an RPI $X_\mathrm{s} = \Poly(F,y_\mathrm{s})$, as above, while the feedback law, too, is parameterized as $\mu_k(x) = u_k + K(x-x_k)$. After pre-computing the support function values
\[
\overline y_i \defeq \support{(A+BK)X_\mathrm{s}}(F_i) \quad \text{and} \quad \overline w_i \defeq \support{\mathbb W}(F_i),
\]
the sequence of homothetic polytopes
\[
X_k = \{ x_k \} + \alpha_k \Poly(F,y_\mathrm{s}),
\]
with $\alpha_k \geq 0$, is a feasible RCT if
\begin{align}
\label{eq::HomotheticRCT}
\hspace{-0.4cm}
\begin{array}{l}
F(Ax_k+Bu_k) + \alpha_k \overline y + \overline w \leq F x_{k+1} + \alpha_{k+1} y_\mathrm{s}
\\
\{ x_k \} + \alpha_k X_\mathrm{s} \subseteq \mathbb X, \ \ \text{and} \ \ \{ u_k \} + \alpha_k K X_\mathrm{s} \subseteq \mathbb U;
\end{array}
\end{align}
see~\cite{Rakovic2012,Rakovic2013} for implementation details and variants. Similar to the rigid RCT approach, the above convex condition is merely sufficient for $X$ to be a feasible RCT because neither the feedback law nor the tube are fully parameterized.

\item \textit{Elastic RCT} parameterizations~\cite{Rakovic2016} are based on the same pre-computations as the above approaches, but, additionally, one precomputes the componentwise minimizer of $\overline Z y_\mathrm{s}$ subject to the dual feasibility constraints
\[
\overline Z F = F(A+BK) \quad \text{and} \quad \overline Z \geq 0 .
\]
Note that this construction is such that, due to weak duality,
\[
(A+BK)\Poly(F,y) \ \subseteq \ \Poly(F, \overline Z y).
\]
This polytopic outer approximation turns out to be tight for all $y = \alpha y_\mathrm{s}$, with $\alpha \geq 0$; see~\cite{Rakovic2016}. In fact, this statement can be further strengthened. Namely, the above outer approximation is also tight for all $y \in \mathbb C(F,y_\mathrm{s})$, as the parametric support function of polytopes is affine on the configuration cone $\mathbb C(F,y_\mathrm{s})$. Unfortunately, however, for general parameters $y$, this overestimation is not tight. Similarly, assuming that $\mathbb X = \Poly(F,\overline y)$ and $\mathbb U = \Poly(G,\overline u)$ are polyhedra with a given facet representation, one can pre-compute the componentwise minimizer of $\widetilde Z y_\mathrm{s}$ subject to 
\[
\widetilde Z F = G K \quad \text{and} \quad \widetilde Z \geq 0,
\]
such that $\widetilde Z y \leq \overline u$ implies $K \Poly(F,y) \subseteq \mathbb U$. The sets $X_k = \{ x_k \} + \Poly(F,y_k)$ are then a feasible RCT if
\begin{align}
\label{eq::ElasticRCT}
\begin{array}{l}
F(Ax_k+Bu_k) + \overline Z y_k + \overline w \leq F x_{k+1} + y_{k+1} \\
F x_k + y_k \leq \overline y \quad \text{and} \quad G u_k + \widetilde Z y_k \leq \overline u.
\end{array}
\end{align}
for all $k \in \mathbb N$; see~\cite{Rakovic2016}. Note that the above condition is convex but merely sufficient for $X$ to be a feasible RCT. This is because the multipliers $\overline Z$ and $\widetilde Z$ are pre-computed. Otherwise, if one regards $\overline Z$ and $\widetilde Z$ as parameters, no conservatism is introduced, but the conditions in~\eqref{eq::ElasticRCT} are bilinear and, hence, non-convex.

\item \textit{Configuration-Constrained RCT} parameterizations~\cite{Villanueva2024} are based on a pre-computation of the support function values
\[
\overline w_i \defeq \support{\mathbb W}(F_i),
\]
and a configuration cone $\Poly(E,0)$, given in facet representation. The configuration-constrained polytopes $\Poly(F,y_k)$ with $E y_k \leq 0$ are a feasible RCT if and only if there exist vertex inputs $U_{k,i}$ such that
\begin{align}
\label{eq::CCRCT}
\begin{array}{l}
F V_i y_k + F B U_{k,i} + \overline w \leq y_{k+1}, \\
U_{k,i} \in \mathbb U, \ V_i y_{k+1} \in \mathbb X, E y \leq 0, \ E y^+ \leq 0
\end{array}
\end{align}
for all $k \in \mathbb N$ and all $i \in \{ 1, \ldots, \nv \}$; see also~\cite{Villanueva2024} for details. This condition is necessary, sufficient, and jointly convex in the parameter sequences $y$ and $U$.

\item \textit{Zonotopic RCT} parameterizations assume that the uncertainty set $\mathbb W = \Zono(G, 0,\overline w)$ is a given zonotope. The sequence of zonotopes
\[
X_k = \Zono(G,c_k,y_k)
\]
is then a feasible RCT if and only if there exist auxiliary variables $z_{k+1,i}$ and vertex controls $U_{k,i}$ with
\begin{align}
\label{eq::ZRCT}
\begin{array}{l}
G z_{k+1,i} + c_{k+1} = B U_{k,i} \\
\phantom{G x_{k+1,i} \ \ } + A \left( G \cdot \mathrm{diag}(y_k) \cdot \sigma_i(G) + c_k \right), \\
\overline w-y_{k+1} \leq z_{k+1,i} \leq y_{k+1}-\overline w, \ U_{k,i} \in \mathbb U, \\
y_k \geq 0, \ G \cdot \mathrm{diag}(y_k) \cdot \sigma_i(G) + c_k \in \mathbb X.
\end{array}
\end{align}
for all $k \in \mathbb N$ and all $i \in \{ 1, \ldots, \nv \}$. This condition is necessary, sufficient, and jointly convex in the parameter sequences $c,y,z$, and $U$.

\end{itemize}
All listed RCT parameterizations restrict either polytopes, feedback laws, or both. Differences between them are discussed in the robust model predictive control literature, as concisely surveyed in Section~\ref{sec::RMPC}.

\subsection{Robust Control Invariant Sets}
\label{sec::RCI}
A set $X$ is called a robust control invariant (RCI) set if $(X,X)$ is an RCT. If $X$ additionally satisfies $X \subseteq \mathbb X$, it is called a feasible RCI set. Similar to RCTs, RCI sets are uniquely characterized by Kerrigan's condition, compare~\eqref{eq::Kerrigan}. Namely, a non-empty set $X$ is an RCI set if and only if there exists a non-empty intermediate set $S \neq \varnothing$ such that
\begin{eqnarray}
\label{eq::KerriganInvariant}
A X \subseteq S + (-B)\mathbb{U}  \qquad \text{and} \qquad S + \mathbb{W}  \subseteq X.
\end{eqnarray}
If $\mathbb X$, $\mathbb U$, and $\mathbb W$ are closed and convex, trivial yet important consequences of Kerrigan's condition can be summarized as follows.

\begin{enumerate}

\item  If $X$ is a feasible RCI set, then its closure is a feasible RCI set.

\item  If $X$ is a feasible RCI set, then its convex hull is a feasible RCI set.

\item If $X$ and $X'$ are feasible RCI sets, then $X \cup X'$ is a feasible RCI set.

\item If $X$ and $X'$ are closed and convex feasible RCI sets, then $\vartheta X + (1-\vartheta)X'$ is a feasible RCI set, $\vartheta \in [0,1]$.

\item If $\mathbb X$ is bounded and $X \neq \varnothing$ a non-empty feasible RCI set, then there exists a closed, convex, and feasible RCI set $X_{\max}$ that contains all other feasible RCI sets. The set $X_{\max}$ is called the maximal RCI set.

\end{enumerate}
While the first four statements are immediate consequences of~\eqref{eq::Kerrigan}, the last statement follows by first using the third statement to show that the union of all feasible RCI sets is a feasible RCI set. The claim follows then from the first and second statements.

Initial endeavors to characterize polytopic RCI sets date back to the second half of the 20th century, with notable reviews conducted by Blanchini in~\cite{Blanchini1999b} and \cite{Blanchini2008}. A significant portion of these early approaches relied on non-convex analysis techniques. This is in contrast to the proliferation of convex analysis-based RCI conditions that have emerged in the past two decades. These convex approaches can be summarized as follows.

\begin{itemize}

\item Conditions~\eqref{eq::RigidRCT},~\eqref{eq::HomotheticRCT},~\eqref{eq::ElasticRCT},~\eqref{eq::CCRCT}, and~\eqref{eq::ZRCT}, here applied to the "short" tube $(X_0,X_1)$, can all be used to derive convex RCI conditions for the respective set parameterization by enforcing stationarity, such that $X_0 = X_1$. This leads to a whole zoo of practical convex conditions for RCIs, each with its respective advantages and disadvantages.

\item Condition~\eqref{eq::SeparableRCT} can be modified in order to derive a convex RCI condition, too. For this aim, however, one needs to first compute an inclusive remainder bound such that
\[
X \ \defeq \ \sum_{j=1}^N \cvx{Y_{N,j}} + \alpha^N \cdot \cvx{Y_\mathrm{s}}
\]
can be proven to be a feasible RCI set. Here, the assumption is that we can find an $\alpha \in [0,1)$ and a pre-computed feasible robust positive invariant polytope $X_\mathrm{s} = \cvx{Y_\mathrm{s}} \subseteq \mathbb X$ in vertex representation, with
\begin{eqnarray}
(A + B K) \cdot X_\mathrm{s} + \mathbb W &\subseteq& X_\mathrm{s} \notag \\
\text{and} \qquad (A + B K) \cdot X_\mathrm{s} &\subseteq& \alpha \cdot X_\mathrm{s} \notag
\end{eqnarray}
for a pre-computed asymptotically stabilizing linear control gain $K$, and pre-computed support values
\begin{align*}
[y_{\mathrm{s}}]_i \ \defeq \ \support{X_\mathrm{s}}( F_i ) \quad
\text{and} \quad [v_{\mathrm{s}}]_i \ \defeq \ \support{K X_\mathrm{s}}( G_i ).
\end{align*}
The above expression for $X$ is then a feasible RCI for $\mathbb X = \Poly(F,\overline y)$, $\mathbb U = \Poly(G,\overline u)$, and $\mathbb W = \cvx{W}$, if there exist slack variables $Y_{k,j}$, $U_{k,j}$, $y_{k,j}$, $v_{k,j}$ with
\begin{align}
\begin{array}{l}
Y_{k+1,j} = A Y_{k,j} + B U_{k,j}, \ \ Y_{j,j} = W, \\
F Y_{k,j} \leq y_{k,j} \one^\tr, \ \ G U_{k,j} \leq v_{k,j} \one^\tr, \\
\alpha^k y_\mathrm{s} + \overset{k}{\underset{j=1}{\sum}} y_{k,j} \leq \overline y, \ \ \alpha^k v_\mathrm{s} + \overset{k}{\underset{j=1}{\sum}} v_{k,j} \leq \overline u
\end{array}
\end{align}
for all $k \in \{ 0, \ldots, N \}$ and all $j \in \{ 1,\ldots, k \}$. Here, we set $Y_{0,0} = 0$, recalling that the limit behavior of an asymptotically stabilized system does not depend on the choice of the initial value. Note that the origins of the above constructions can (in a similar version) be traced back to the seminal work of Rakovi\'c, Kerrigan, Mayne, and many other co-workers; see~\cite{Rakovic2007} and later also~\cite{Rakovic2010},~\cite{Rakovic2011FP}, and~\cite{Rakovic2012a}, where various implementation details, variants, and derivations can be found.
\end{itemize}

In summary, Conditions \eqref{eq::SeparableRCT}--\eqref{eq::ZRCT} and their variants provide a solid foundation for synthesizing convex RCT and RCI parameterizations. While each method has unique strengths and limitations, none is perfect, as they all rely on parameterizing sets, feedback laws, or both. For instance, affine feedback sub-optimality is analyzed in~\cite{Bertsimas2012}. For a broader perspective on RCTs in nonlinear viability theory, see~\cite{Aubin2009}.

\subsection{Control Lyapunov Functions}
\label{sec::CLF}
Control Lyapunov Functions (CLFs)~\cite{Zubov1965} are at the core of nonlinear optimal control theory~\cite{Artstein1983,Bellman1957} and dynamic programming~\cite{Bertsekas2012}, rendering them a broad and comprehensive topic that exceeds the scope of this review. Nevertheless, there are certain conceptual ideas that are relevant in the context of polyhedral control. These definitions focus on the special case that $\mathbb X$ and $\mathbb U$ are closed and convex, $(0,0) \in \mathbb X \times \mathbb U$, and $L: \mathbb R^{n_x} \times \mathbb R^{n_u} \to \mathbb R$ is a convex, non-negative stage cost with $L(0,0) = 0$.

\begin{itemize}

\item A function $M: \mathbb R^n \to \mathbb R \cup \{ \infty \}$, is called an $L$-CLF if it is non-negative, proper compact convex, $M(0) = 0$ and $M(x) = \infty$ for all $x \notin \mathbb X$, and
\begin{align}
\notag
M(x) \ \geq \ \min_{u \in \mathbb U} \ \left\{ L(x,u) + M(Ax+Bu) \right\}.
\end{align}

\item An $L$-CLF is called piecewise affine if its epigraph is a polyhedron; see~\cite{MolchanovI,MolchanovII,MolchanovIII} as well as~\cite{Nguyen2018}.

\end{itemize}
The domain of an $L$-CLF function $M$, which is defined as $D(M) \defeq \{ x \mid M(x) < \infty \}$, and all of the sub-level set the $M$ are, by construction, closed and convex control invariant sets, see also~\cite{Artstein1983,Blanchini2008} or~\cite{Giesl2015} for a formal proof. Consequently, one can pick at least one associated feasible control law,
\[
\mu[M,L](x) \ \in \ \underset{u \in \mathbb U}{\text{argmin}} \ \left\{ L(x,u) + M(Ax+Bu) \right\},
\]
which is then well-defined for all $x$ in the iterated domain of $M$. This iterated domain is denoted by
\[
D^+(M) \ \defeq \ \left\{ \ x \in \mathbb X \ \middle|
\begin{array}{l}
\exists u \in \mathbb U: \\
M(Ax+Bu) < \infty
\end{array}
\right\} \ \supseteq \ D(M).
\]
The latter inclusion follows from the fact $D(M)$ is control invariant and, consequently, $D^+(M)$ is control invariant. The associated closed-loop system
\[
x_{k+1} = A x_k + B u_k \quad \text{with} \quad u_k = \mu[M,L](x_k)
\]
satisfies $x_k \in \mathbb X$ as well as $u_k \in \mathbb U$ for all $k \in \mathbb N$ as well as
\[
\sum_{k=0}^\infty L(x_k,u_k) \ \leq \ M(x_0)
\]
for all $x_0$ with $M(x_0) < \infty$. This statement follows under the above mentioned assumptions on $\mathbb X$, $\mathbb U$, and $L$ by applying a simple induction argument; see also~\cite{Giesl2015} for details. Moreover, the above infinite horizon performance bound can be used as a starting point to deduce various statements about the stability or convergence behavior of the closed-loop system under $\mu[M,L]$. For example, since $L$ is convex and non-negative, the zero-level set
\[
X_L \ \defeq \ \left\{ \ x \ \middle| \ 0 \ = \ \min_{u \in \mathbb U} L(x,u) \ \right\}
\]
is closed and convex. Hence, since convex functions are Lipschitz continuous, it follows that the trajectories of the closed-loop system converge to $X_L$,
\[
\lim_{k \to \infty} \ \min_{x \in X_L} \| x_k - x \| \ = \ 0,
\]
and it also follows from the above definitions that the $L$-CLF $M$ necessarily satisfies $M(x) > 0$ for all $x \notin X_L$. In particular, for the special case that $X_L = \{ 0 \}$ is a singleton, $M$ must be positive definite, and the associated closed-loop system is asymptotically stable \cite{Zubov1965}. Moreover, under the above mentioned assumptions on $\mathbb X$, $\mathbb U$, and $L$, the following statements hold.

\begin{enumerate}

\item A closed and convex set $X$ with $0 \in X$ is control invariant if and only if its (modified) gauge function,
\[
\mathcal G(x) \ \defeq \ \min_{\alpha \in [0,1]} \alpha \ \ \text{s.t.} \ \ x \in \alpha X,
\]
is a $0$-CLF.

\item If $M_k$ is an $L$-CLF, then the function
\[
M_{k+1}(x) \ \defeq \ \min_{u \in \mathbb U} \ \left\{ L(x,u) + M_k(Ax+Bu) \right\}
\]
is also an $L$-CLF. It satisfies $M_{k+1} \leq M_k$. Moreover, we have $D(M_{k+1}) = D^+(M_k)$.

\item If there exists at least one $L$-CLF for a given system, then there exists a unique $L$-CLF, denoted by $M_\infty$, whose epigraph contains the epigraph of all other robust $L$-CLFs. It satisfies the stationary discrete-time Hamilton-Jacobi-Bellman (HJB) equation
\[
M_\infty(x) \ = \ \min_{u \in \mathbb U} \ \left\{ L(x,u) + M_\infty(Ax+Bu) \right\}.
\]
Furthermore, the set $D^+(M_\infty) = D(M_\infty)$ is control invariant. It is the maximal control invariant set, as any control invariant set is a subset of $D(M_\infty)$.
\end{enumerate}
The gauge function $\mathcal G$, also known as the "Minkowski-Lyapunov function" associated with the control invariant set $X$ or simply "Minkowski CLF"~\cite[Thm.4.24]{Blanchini2008}, builds upon pioneering work in~\cite{Kiendl1992} and~\cite{Sznaier1993} exploring norms as Lyapunov functions. For a broader understanding of Minkowski-CLF properties, we refer to~\cite{Rakovic2020}, while practical piecewise affine constructions are detailed in~\cite{Rakovic2024}. The monotonous CLF recursion in Item~2 above, akin to the discrete-time dynamic programming recursion for linear systems~\cite{Bertsekas2012}, generates a monotonically decreasing sequence of $L$-CLFs $M_0,M_1,M_2,\ldots$ that converge to an $M_\infty$ satisfying the discrete-time HJB equation. This convergence analysis leverages standard convex analysis techniques from~\cite{Rockafellar1970} and is elaborated in~\cite{Bertsekas2012}.

A direct consequence of the Minkowski CLF construction is that all reviewed methods for computing polyhedral control invariant sets, including those for RCI sets, can readily yield $0$-CLF functions. Early polyhedral computing methods based on similar constructions of piecewise affine Minkowski CLFs were introduced by Blanchini and colleagues in the 1990s~\cite{Blanchini1994,Blanchini1995,Blanchini1999a}. More recent advancements, exemplified in~\cite{Rakovic2024}, have extended these principles to construct piecewise affine $L$-CLFs, often tailored to $L$ being a Minkowski gauge function. For a comprehensive overview of iterative and approximate dynamic programming techniques, we refer to~\cite{Giesl2015}.

Finally, in terms of convex optimization-based non-iterative methods that directly compute a piecewise affine CLF approximation of the solution $M_\infty$ of the stationary Hamilton-Jacobi equation, one option is to work with configuration-constrained epigraph representations. Here, the main idea is to work with polyhedra $\Poly(F,y) \subseteq \mathbb R^n$, with $n = n_x+1$, that are parameterized over a suitable configuration cone $\Poly(E,0)$ using vertex matrices $V_1,\ldots,V_{\nv}$, such that
\begin{align}
\label{eq::PolyCLF}
\Poly(F,y) \ = \ \mathrm{conv} \left( V_1 y, \ldots, V_{\nv} y \right) + \cone{[0^\tr,1]^\tr}
\end{align}
is for all $y \in \Poly(E,0)$ a Minkowski-Weyl representation of the epigraph of an associated piecewise-affine function on a polytopic domain. In this context, we may assume, without loss of generality, that the vertices are enumerated in such a way that the normal cone of the vertex $V_1 y$ contains the downward pointing vector $[0^\tr,-1]^\tr$. If we use such an enumeration, $\Poly(F,y)$ with $E y \leq 0$ is the epigraph of an $L$-CLF if and only if there exist vertex control inputs $U_i$ and slack variables $\omega_i$ such that
\begin{align}
\label{eq::CCCLF}
\hspace{-0.4cm}
\begin{array}{l}
L( \Omega V_i y, U_i ) + \omega_i \leq y_n, \ \ \Omega V_i y \in \mathbb X, \ \ U_i \in \mathbb U, \\[0.1cm]
F \left(
\begin{array}{c}
A \Omega V_i y + B U_i \\
\omega_i
\end{array}
\right)
 \leq y, \ V_1 y = 0, \ \ \text{and} \ \ E y \leq 0
\end{array}
\end{align}
for all $i \in \{ 1, \ldots, \nv \}$, where $\Omega$ denotes the first $n_x$ rows of the $n \times n$ unit matrix.

\subsection{Robust Control Lyapunov Functions}
\label{sec::RCLF}
Let $\mathbb W$ be a compact set, $\mathbb X$ and $\mathbb U$ closed and convex, and let $X_\mathrm{s}$ be a feasible RCI set. In such a setting, a non-negative and convex cost function $L: \mathbb R^{n_x} \times \mathbb R^{n_u} \to \mathbb R$ is an admissible tracking cost for $X_\mathrm{s}$ if
\[
\forall x \in X_\mathrm{s}, \ \exists u \in \mathbb U: \ \forall w \in \mathbb W, \quad
\left\{
\begin{array}{l}
L(x,u) = 0 \\
Ax + B u + w \in X_\mathrm{s}.
\end{array}
\right.
\]
This condition necessarily implies that the zero-level set of $L$ satisfies
\[
X_L \ = \ \left\{ \ x \in \mathbb X \ \middle| \ \min_{u \in \mathbb U}  L(x,u) \ = \ 0 \ \right\} \ \supseteq \ X_\mathrm{s}.
\]
Additionally, if $X_L = X_\mathrm{s}$, we say that $L$ is a positive definite stage cost function with respect to the target set $X_\mathrm{s}$. Here, in principle, any of the above reviewed polyhedral computing methods for constructing RCI sets can be used to construct the target set $X_\mathrm{s}$, see also~\cite[Sect.~3.2]{Kerrigan2004} for a discussion of this aspect in the context of min-max optimal control for linear systems.

Next, in analogy to the nominal case non-negative function $M: \mathbb R^n \to \mathbb R \cup \{ \infty \}$, is called a robust $L$-CLF if it is non-negative, proper compact convex, $M(x) = 0$ for all $x \in X_\mathrm{s}$ and $M(x) = \infty$ for all $x \notin \mathbb X$, and
\begin{align}
\label{eq::robustCLF}
M(x) \ \geq \ \min_{u \in \mathbb U} \max_{w \in \mathbb W} \ \left\{ L(x,u) + M(Ax+Bu+w) \right\}.
\end{align}
Apart from the above apparent complications regarding the construction of an admissible tracking cost $L$, the analysis of robust $L$-CLFs is essentially analogous to the analysis of their nominal counterparts. In particular, the sublevel sets of an $L$-CLF $M$, including its domain $D(M)$, are convex and closed RCIs. Moreover, on the domain
\[
D^+(M) \ \defeq \ \left\{ \ x \in \mathbb X \ \middle|
\begin{array}{l}
\exists u \in \mathbb U: \ \forall w \in \mathbb W \\
M(Ax+Bu+w) < \infty
\end{array}
\right\},
\]
a parametric minimizer $\mu[M,L]$ of the right-hand of~\eqref{eq::robustCLF} exists. It then follows again from a simple induction argument~\cite{Giesl2015} that the closed-loop trajectories,
\[
x_{k+1} \ = \  Ax_k + B u_k + w_k \quad \text{with} \quad u_k = \mu[M,L](x_k) 
\]
satisfy
\[
\sum_{k=0}^\infty L(x_k,u_k) \ \leq \ M(x_0)
\]
for all $x_0 \in D^+(M)$, independent of the uncertainty realization $w_0,w_1,\ldots \in \mathbb W$. Additionally, the min-max dynamic programming recursion
\[
M_{k+1}(x) \ \defeq \ \min_{u \in \mathbb U} \max_{w \in \mathbb W} \ \left\{ L(x,u) + M_k(Ax+Bu+w) \right\}
\]
can be used to generate a sequence of robust $L$-CLFs, starting from an initial $L$-CLF $M_0$. Hence, under the above mentioned assumptions, if such an initial $L$-CLF exists, the limit $M_\infty$ of the above sequence exists, is itself and $L$-CLF, and satisfies the discrete-time min-max Hamilton-Jacobi-Bellman equation
\[
M_\infty(x) \ = \ \min_{u \in \mathbb U} \max_{w \in \mathbb W} \ \left\{ L(x,u) + M_\infty(Ax+Bu+w) \right\}.
\]
Note that min-max dynamic programming approaches, rooted in Witsenhausen's seminal work from the late 1960s \cite{Witsenhausen1968}, have long been established. However, methods leveraging these approaches for approximating $M_\infty$ emerged significantly later. Early iterative polyhedral techniques, maintaining bounds on $M_\infty$, were presented in \cite{Diehl2004} for piecewise affine CLF approximations. Alternatively, approximations with pre-defined complexity can be obtained via a single convex optimization problem. For example, by using a similar construction as previously, the configuration-constrained polyhedron in \eqref{eq::PolyCLF} serves as the epigraph of a robust CLF function on $\Poly(E,0)$ if and only if
\begin{align}
\notag
\hspace{-0.4cm}
\begin{array}{l}
L( \Omega V_i y, U_i ) + \omega_i \leq y_n, \ \ \Omega V_i y \in \mathbb X, \ \ U_i \in \mathbb U, \\[0.1cm]
F \left(
\begin{array}{c}
A \Omega V_i y + B U_i \\
\omega_i
\end{array}
\right) + \overline w
 \leq y, \ y_1 = 0, \ \ \text{and} \ \ E y \leq 0
\end{array}
\end{align}
for all $i \in \{ 1, \ldots, \nv \}$. Here, we assume that $F_1 = [0^\tr,-1]^\tr$ is a downward pointing facet normal and
\[
\overline w_i \ = \ \support{ \{ (w,0)^\tr \mid w \in \mathbb W \} }(F_i)
\]
pre-computed support function values.

\subsection{Model Predictive Control}
\label{sec::MPC}
Instead of explicitly pre-computing CLFs that approximate the solution of HJBs with high precision, an alternative approach is to start with an arbitrary $L$-CLF $M_0$, but actually use the implicit representation
\begin{align}
\label{eq::MPC}
M_{N}(\hat x) \ = \ & \min_{x,u} \ \sum_{k=0}^{N-1} L(x_k,u_k) + 
M_0(x_N) \\
 & \ \mathrm{s.t.} \ \left\{
\begin{array}{l}
\forall k \in \{ 0, \ldots, N-1 \} \\
x_{k+1} = A x_k + B u_k, \\
x_k \in \mathbb X, \ u_k \in \mathbb U, \ x_0 = \hat x.
\end{array}
\right. \notag
\end{align}
for the function $M_N$ that is obtained by propagating $M_0$ via the dynamic programming recursion, $N$ steps backward in time. In this context, the integer $N \geq 1$ is called the prediction horizon. Here, we recall that if $M_0$ is an $L$-CLF, then $M_{N}$ is an $L$-CLF too, but $M_N$ is a better approximation of the infinite horizon cost $M_\infty$. This is because we have
\[
0 \ \leq \ L \ \leq \ M_\infty \ \leq \ M_{N}  \ \leq \ M_0.
\]
Controllers that use such implicit CLF representations have become popular under the umbrella name "\textit{model predictive control}" (MPC)~\cite{Kouvaritakis2015} and~\cite{Rawlings2015}. They are typically implemented by solving the optimization problem~\eqref{eq::MPC} online to compute the parametric minimizer, $(x^\star(\hat x),u^\star(\hat x))$, directly upon receiving the state measurement $\hat x$. This means that the actual MPC control law is given by
\[
\mu[M_{N-1},L](\hat x) \ = \ u_0^\star(\hat x),
\]
recalling that this control law is well-defined on the control invariant domain $D^+(M_{N-1}) = D(M_N)$; see also~\cite{Mayne2000}. Two simple yet important consequences of this construction can be formulated as follows~\cite{Rawlings2015}.

\begin{enumerate}

\item Because $D^+(M_{N-1}) = D(M_N)$ is control invariant, the MPC closed-loop trajectories
\begin{align}
\label{eq::MPCtrajectories}
x_{k+1} = A x_k + \mu[M_{N-1},L](x_k),
\end{align}
started at $x_0 \in D(M_N)$, satisfy $x_k \in D(M_N)$ for all $k \in \mathbb N$. In other words, the MPC controller~\eqref{eq::MPC} is recursively feasible.

\item The MPC closed-loop trajectories, as generated by~\eqref{eq::MPCtrajectories}, converge into the zero level set of $L$,
\[
\lim_{k \to \infty} \ \min_{x \in X_L} \| x_k - x \|,
\]
as long as $x_0 \in D(M_N)$. 
\end{enumerate}
Note that many different MPC frameworks and variants of~\eqref{eq::MPC} have been proposed in the last three decades, such as using additional terminal constraints together with an $L$-CLF $M_0$ that is only defined locally, or completely omitting $M_0$ in the formulation of~\eqref{eq::MPC}. A more detailed introduction and overview of this large research field can be found in the textbooks~\cite{Kouvaritakis2015,Rawlings2015,Gruene2016} or the survey article~\cite{Mayne2014,Kohler2024}. Furthermore, we note that extensions to more general economic stage cost functions have been developed under suitable dissipativity assumptions, see, for example, ~\cite{Faulwasser2018,Houska2017}.

\subsection{Robust Model Predictive Control}
\label{sec::RMPC}
Implicit representations for CLFs can also be used to approximate the solution of min-max HJBs. Similar to the nominal case, these methods start with a robust $L$-CLF $M_0$ and write out the associated $N$-step backward dynamic programming recursion by enumerating all possible scenarios~\cite{Kerrigan2004}. By using the same notation as in Condition~\eqref{eq::ScenarioRCT}, the function
\begin{eqnarray}
\label{eq::ScenarioMPC}
\hspace{-0.4cm}
M_N( \hat x) &=& \min_{x,u} \ \max_{i \in \mathcal I_N} \ \sum_{k=0}^{N-1} L(x_k^i,u_{s(k,i)}) + M_0(x_N^i) \notag \\
& & \ \mathrm{s.t.} \ \left\{
\begin{array}{l}
\forall k \in \{0,\ldots N-1\}, \ \forall i \in \mathcal I_N, \\
x_{k+1}^i = A x_k^i + B u_{s(k,i)} + w_k^i \\
x_k^i \in \mathbb X, \ u_{s(k,i)} \in \mathbb U, \
x_0^i = \hat x
\end{array}
\right.
\end{eqnarray}
is an $L$-CLF that approximates $M_\infty$ more accurately than $M_0$ does, as $0 \leq L \leq M_\infty \leq M_N \leq M_0$. Controllers that evaluate the associated control laws $\mu[M_{N-1},L]$ by solving~\eqref{eq::ScenarioMPC} online are called "\textit{Worst-Case Scenario MPC}" controllers or simply "\textit{Min-Max MPC}" controllers; see also~\cite{Bemporad2003},~\cite{Diehl2004},~\cite{Lofberg2003},~\cite{Lofberg2003a},  and~\cite{Scokaert1998} for variants. Note that such Min-Max MPC controllers \mbox{are---by} \mbox{construction---recursively} feasible on the domain $D(M_N)$, recalling that the sub-level sets of $M_N$, including its domain, are RCI sets. Moreover, they are robustly convergent in the sense that all trajectories of the closed-loop system
\[
x_{k+1} \ = \ A x_k + B \mu[M_{N-1},L](x_k) + w_k 
\]
converge, independent of the uncertainty realization, to $X_L$, the $0$-level set of $L$, recalling that this set coincides with the target RCI set $X_\mathrm{s}$ whenever $L$ is positive definite with respect to this target.

A survey on related stability properties of min-max MPC controllers is presented in~\cite{Raimondo2009}. Additionally, Mayne critiques robust MPC in~\cite{Mayne2015}, suggesting nominal MPC's robustness, as studied in~\cite{Pannocchia2011}, may suffice for many applications. In scenario-based MPC, this critique is partially valid due to the exponential growth of scenarios with prediction horizon, limiting the practicality of such scenario-based implementations. However, as we will survey below, various efficient polynomial run-time formulations of robust MPC have emerged, partially addressing Mayne's concerns.

In recent years, Tube MPC controllers have gained popularity as a viable alternative class of robust MPC controllers. These methods are rooted in the concept of propagating RCTs in the state space, first introduced in the 1970s, as exemplified in~\cite{Bertsekas1972}; see also~\cite{Kothare1996} and~\cite{Cuzzola2002}. Modern formulations of Tube MPC typically entail solving optimization problems of the form
\begin{align}
\label{eq::TMPC}
J_N(\hat x) \ \defeq & \min_{X \in \mathcal X,\mu \in \mathcal U} \ \sum_{k=0}^{N-1} \mathcal L( X_k, \mu_k(X_k) ) + \mathcal M_0(X_N) \\
& \ \ \ \text{s.t.} \ \ \ \left\{
\begin{array}{l}
\forall k \in \{ 0, \ldots, N-1 \}, \\
\forall x_k \in X_k, \forall w_k \in \mathbb W, \\ \ Ax_k+B \mu_k(x) + w_k \in X_{k+1}, \\
X_k \subseteq \mathbb X, \ \mu_k(X_k) \subseteq \mathbb U, \ \hat x \in X_0.
\end{array}
\right. \notag 
\end{align}
In this context, $\mathcal X$ denotes a suitable set of RCTs, and $\mathcal U$ denotes a suitable set of ancillary control law sequences, where
\[
\mu_k(X_k) \ \defeq \{ \ \mu_k(x_k) \ \mid \ x_k \in X_k \ \}
\]
denotes the elements of the control tube. In order to discuss conditions under which the above general tube MPC controller is well-formulated, we first recall that a function $\mathcal M: 2^{\mathbb R^n} \to \mathbb R \cup \{ \infty \}$ is called
\begin{itemize}
\item monotonous if $X \subseteq X'$ implies $\mathcal M(X) \leq \mathcal M(X')$,
\item proper if there exists an $X \neq \varnothing$ with $M(X) < \infty$,
\item closed if $\mathcal M(X) = \mathcal M( \mathrm{cl}(X) )$ for all sets $X$,
\item and convex if both $\mathcal M(X) = \mathcal M( \cvx{X} )$ and
\[
M(\vartheta X + (1-\vartheta) X') \leq \vartheta M(X) + (1-\vartheta) M(X')
\]
for all $X,X' \subseteq \mathbb R^n$ and all $\vartheta \in [0,1]$.
\end{itemize}
Assuming monotonicity, properness, closedness, and joint convexity in set-valued variables for
\[
\mathcal L: 2^{\mathbb R^{n_x}} \times 2^{\mathbb R^{n_u}} \to \mathbb R \cup \{ \infty \} \quad \text{and} \quad \mathcal M_0: 2^{\mathbb R^{n_x}} \to \mathbb R \cup \{ \infty \},
\]
Kerrigan's condition~\eqref{eq::Kerrigan} for RCTs, for closed and convex $\mathbb X$, $\mathbb U$, and $\mathbb W$, allows restricting $X$ to element-wise closed and convex sets without loss of generality. This stems from feasibility preservation and objective value equivalence after taking element-wise convex hulls. To illustrate the framework's modeling power, consider applications such as set diameter computation, least-squares tracking error against a target set $X_\mathrm{s}$, or evaluating the worst-case value of proper, closed, and convex stage functions $\ell$,
\begin{eqnarray}
\mathcal M(X) &=& \sup_{x,x' \in X} \| x-x'\| \notag \\[0.16cm]
\mathcal M(X) &=& \sup_{x \in X} \ \underset{x' \in X_\mathrm{s}}{\inf\vphantom{p}} \ \| x-x' \|^2 \notag \\[0.16cm]
\text{or} \qquad \mathcal M(X) &=& \sup_{x \in X} \ \ell(x), \notag
\end{eqnarray}
which are all examples of monotonous, proper, closed, and convex functions. These, along with numerous other classes of geometrically meaningful functions, can be used for designing the stage cost $\mathcal L$. However, when selecting the terminal cost $\mathcal M_0$, additional considerations pertaining to stability become essential. Assuming that $L$ is a non-negative, monotonous, proper, closed, and convex stage cost, a function $\mathcal M$ that is also non-negative, monotonous, proper, closed, and convex is called an $\mathcal L$-CLF function if
\begin{align}
\mathcal M(X) \ \geq \ & \min_{\mu,X^+} \ \mathcal L(X,\mu(X)) + \mathcal M(X^+) \notag \\
& \ \text{s.t.} \ \ \left\{
\begin{array}{l}
\forall x \in X, \forall w \in \mathbb W, \\ \ Ax+B \mu(x) + w \in X^+, \\
X,X^+ \subseteq \mathbb X, \ \mu_k(X) \subseteq \mathbb U .
\end{array}
\right.
\end{align}
Once this definition is introduced, the subsequent discussion of the properties of Tube MPC problems becomes largely analogous to that of nominal MPC. The primary differences between the tube MPC problem~\eqref{eq::TMPC} and the nominal problem~\eqref{eq::MPC} lie in the replacement of the terminal cost $M_0$ in~\eqref{eq::MPC} with the $\mathcal L$-CLF function $\mathcal M_0$, the replacement of $L$ with $\mathcal L$, and the substitution of states and controls, $x_k$ and $u_k$, with sets $X_k$ and ancillary control laws $\mu_k$, respectively. Specifically, if $\mathcal M_0$ is an $\mathcal L$-CLF, then the iterates of the set-theoretic dynamic programming recursion
\begin{align}
\mathcal M_{k+1}(X) \ = \ & \min_{\mu,X^+} \ \mathcal L(X,\mu(X)) + \mathcal M_k(X^+) \notag \\
& \ \text{s.t.} \ \ \left\{
\begin{array}{l}
\forall x \in X, \forall w \in \mathbb W, \\ \ Ax+B \mu(x) + w \in X^+, \\
X,X^+ \subseteq \mathbb X, \ \mu_k(X) \subseteq \mathbb U
\end{array}
\right.
\end{align}
are, by construction, $\mathcal L$-CLF functions, too. For the theoretically relevant case that $\mathcal X$ and $\mathcal U$ contain the set of all possible sequences of closed convex sets and all possible sequences of feedback laws, the value function $J_N$ of~\eqref{eq::TMPC} can be related to $\mathcal M_N$, as the required monotonicity property of $\mathcal L$-CLF function implies that we have
\[
J_N(\hat x) = \mathcal M_N(\{\hat x\}).
\]
In particular, recursive feasibility of~\eqref{eq::TMPC} on the domain of $\mathcal M_N$ and convergence of the closed-loop trajectories to the zero level set
\[
X_{\mathcal L} \ \defeq \ \left\{ \ x \in \mathbb X \ \middle| \ 0 \ = \ \min_{u \in \mathbb U} \ \mathcal L \left( \{x\}, \{u\} \right) \ \right\}
\]
for all possible uncertainty scenarios are then immediate consequences of this construction and our assumption that $\mathcal L$ is monotonous, too. An in-depth discussion of related set-theoretic stability concepts and proofs of asymptotic convergence of Tube MPC for stage cost functions $\mathcal L$ that are positive definite with respect to a given or optimal RCI target set can be found in~\cite{Villanueva2020}.

When parameterizing \eqref{eq::TMPC} with practical polyhedra, any RCT parameterization from Section~\ref{sec::RCT} can transform the problem into a tractable convex optimization problem for specific state- and control tubes $\mathcal X$ and $\mathcal U$. However, careful design of the stage cost $\mathcal L$, and associated $\mathcal L$-CLF function $\mathcal M_0$ is crucial, as they must be compatible with the chosen RCT and feedback law parameterizations. This necessitates tailored terminal cost conditions for each parameterization, which are not inherently encoded in the general $\mathcal L$-CLF definition. Examples of convex optimization-based polyhedral Tube MPC controllers can be listed as follows.

\begin{enumerate}

\item \textit{Homothetic Tube MPC} constructs $\mathcal X$ and $\mathcal U$ by parameterizing the state and control tubes such that
\begin{eqnarray}
X_k &=& \{ x_k \} + \alpha_k \mathcal P( F, y_\mathrm{s} ) \\
\text{and} \quad \mu_k(X_k) &=& \{ u_k \} + \alpha_k K \mathcal P( F, y_\mathrm{s} )
\end{eqnarray}
for parameter sequences $(u_k,x_k,\alpha_k)$, with $\alpha_k \geq 0$, using the same assumptions and notation that were introduced in the context of Condition~\eqref{eq::HomotheticRCT}. This means that, under these assumptions, the constraints in~\eqref{eq::TMPC} are satisfied if Condition~\eqref{eq::HomotheticRCT} holds and
\begin{align}
\label{eq::initialHTMPC}
F (\hat x - x_0) \ \leq \ \alpha_0 y_\mathrm{s}.
\end{align}
Moreover, if $\mathcal L$ is a non-negative, monotonous, proper, closed, and convex stage cost, the function
\begin{align}
\label{eq::LHTMPC}
& L(x_k,\alpha_k,u_k) \\
& \ \defeq \mathcal L( \{ x_k \} + \alpha_k \mathcal P( F, y_\mathrm{s} ), \{ u_k \} + \alpha_k K \mathcal P( F, y_\mathrm{s} ) ) \notag
\end{align}
is jointly convex in $(u_k,x_k,\alpha_k)$, for $\alpha_k \geq 0$, and can be found explicitly for many practical choices of $\mathcal L$~\cite{Rakovic2012}. Hence, an associated parameterized version of~\eqref{eq::TMPC} is given by
\begin{align}
& \min_{(x,\alpha,u)} \ \sum_{k=0}^{N-1} L(x_k,\alpha_k,u_k) + M_0(x_N,\alpha_N) \ \notag \\
& \ \ \ \text{s.t.} \ \ \ 
\text{Conditions~\eqref{eq::HomotheticRCT}~and~\eqref{eq::initialHTMPC}.} \notag
\end{align}
As mentioned above, however, in the context of set parameterizations and their associated, potentially conservative RCT conditions, one needs to adapt the definition $\mathcal L$-CLFs to the particular recursivity conditions. For example, in the particular case of the above homothetic RCT parameterization, $M_0$ is a compliant $L$-CLF if
\begin{align}
\label{eq::HTMPCTerminal}
M_0(x,\alpha) \ \leq \ & \min_{u,\alpha^+,x^+} L(x,\alpha,u) + M_0(x^+,\alpha^+) \\
& \ \text{s.t.} \
\left\{
\begin{array}{l}
F(Ax+Bu) + \alpha \overline y + \overline w \\
\ \leq \ F x^+ + \alpha^+ y_\mathrm{s}
\\
\{ x_k \} + \alpha_k X_\mathrm{s} \subseteq \mathbb X, \\
\{ u_k \} + \alpha_k K X_\mathrm{s} \subseteq \mathbb U.
\end{array}
\right. \notag
\end{align}
If one uses a terminal cost $M_0$ that satisfies this condition, the Homothetic Tube MPC implementation is recursively feasible and its associated closed-loop states converge to the corresponding zero-level set of $L$. The proof of this statement is then, of course, completely analogous to the above statement for the general set-theoretic Tube MPC controller with the only difference being that our definition of $\mathcal L$-CLFs is adapted to the chosen set parameterization.

Note that the above homothetic Tube MPC parameterization can be traced back to the work of Langson and co-workers~\cite{Langson2004}, although a complete analysis and implementation of such a controller has only been worked out later, by Rakovi\'c and co-workers~\cite{Rakovic2012,Rakovic2013}. 

\item \textit{Rigid Tube MPC} can be interpreted as a restricted variant of Homothetic Tube MPC that sets the scaling parameter to $\alpha_k = 1$. The state and control tubes have the form
\begin{eqnarray}
X_k &=& \{ x_k \} + \mathcal P( F, y_\mathrm{s} ) \\
\text{and} \quad \mu_k(X_k) &=& \{ u_k \} + K \mathcal P( F, y_\mathrm{s}  ),
\end{eqnarray}
where Condition~\eqref{eq::RigidRCT} together with~\eqref{eq::initialHTMPC} is used to replace the constraints in~\eqref{eq::TMPC}. Such rigid tube controllers are traditionally implemented by pre-computing constant constraint margins, as in the original article by Mayne and co-workers~\cite{Mayne2005}. Similarly,~\cite{Chisci2001} uses a pre-computed, but time-varying, tube based on reachable sets and corresponding pre-computed constraint margins. Recent implicit implementation of rigid Tube MPC, however, side-track most set-based offline computations by using an implicit constraint margin representation, which has certain advantages upon scaling the implementation of these controllers for problems in higher dimensional spaces~\cite{Rakovic2023}.

\item \textit{Elastic Tube MPC} constructs $\mathcal X$ and $\mathcal U$ by parameterizing the state- and control tubes such that
\begin{eqnarray}
X_k &=& \{ x_k \} + \mathcal P( F, y_k ) \\
\text{and} \quad \mu_k(X_k) &=& \{ u_k \} + K \mathcal P( F, y_k )
\end{eqnarray}
for parameter sequences $(u_k,x_k,y_k)$, with $y_k \geq 0$. The constraints in~\eqref{eq::TMPC} are then parameterized by replacing them with the conservative but convex recursivity condition~\eqref{eq::ElasticRCT} and initial value constraint,
\[
F(\hat x - x_0) \ \leq \ y_0.
\]
Details on the construction of the associated stage and terminal CLF functions can be found in~\cite{Rakovic2016}. Note that, in the context of robust control, the above parameterization has originally been proposed by Cannon and co-workers~\cite{Cannon2015}. The name "Elastic Tube MPC" has then been introduced in a later article by Rakovi\'c and co-workers~\cite{Rakovic2016} in the context of Tube MPC.

\item \textit{Configuration Constrained Tube MPC}, as originally proposed in~\cite{Villanueva2024}, constructs $\mathcal X$ and $\mathcal U$ by parameterizing the state- and control tubes such that
\begin{eqnarray}
X_k &=& \mathcal P( F, y_k ) \\
\text{and} \quad \mu_k(X_k) &=& \cvx{U_k}
\end{eqnarray}
for parameter sequences $(y_k,U_k)$, with $E y_k \leq 0$, using the same assumptions and notation that were introduced in the context of Condition~\eqref{eq::CCRCT}. Details on the construction of stage costs can be found in~\cite{Villanueva2024} and~\cite{Badalamenti2024}, while an in-depth discussion of terminal cost functions for Configuration Constrained Tube MPC can be found in~\cite{Houska2024}.

\item \textit{Separable Feedback MPC}, in one of its most general variants~\cite{Rakovic2011FP,Rakovic2012b}, parameterizes state- and control tubes as
\begin{eqnarray}
\begin{array}{rcl}
X_k &=& \sum_{j=0}^{k} \mathrm{conv}(Y_{k,j}) \\
\text{and} \quad \mu_k(X_k) &=& \sum_{j=0}^{k} \mathrm{conv}(U_{k,j})
\end{array}
\end{eqnarray}
using Condition~\eqref{eq::SeparableRCT} to parameterize the RCT constraints in~\eqref{eq::TMPC}. If one wishes to avoid the explicit constructions of the above Minkowski sums, however, special care needs to be taken upon designing appropriate stage and terminal cost functions, as discussed in~\cite{Rakovic2011FP}. Note that the above general separable feedback scheme is also known under the name~\textit{Fully Parameterized Tube MPC}. Early variants of such Separable Feedback MPC schemes have been developed by using more restrictive disturbance affine feedback parameterization~\cite{Goulart2006}, with original ideas tracing back to~\cite{vanHessem2003}. Recent extensions can be found in~\cite{Buschermohle2024,Nair2022}. Efficient implementations of disturbance affine feedback MPC have also appeared recently under the name \text{System-Level Control}~\cite{Sieber2022}.
\end{enumerate}
Note that numerical comparisons of the aforementioned Tube MPC schemes can be found in~\cite{Villanueva2024}, while extensions of robust MPC for non-tracking objectives are discussed in~\cite{Bayer2014,Broomhead2015,Angeli2021,Kloppelt2021,Schwenkel2020,Villanueva2020}. Additionally, an overview of extensions of Tube MPC to nonlinear and continuous-time systems is provided in~\cite{Villanueva2017}, with a more general survey of nonlinear robust MPC presented in~\cite{Houska2019}.

Furthermore, set-theoretic MPC methods have been extended to output-feedback MPC formulations, as exemplified in \cite{Bemporad2000}, \cite{Brunner2018},~\cite{Efimov2022}, and \cite{Mayne2006}. Drawing from earlier work by Witsenhausen \cite{Witsenhausen1968a}, more general set-theoretic notions have also been developed in \cite{Dorea2009, Dorea2021} for single invariance and in \cite{Artstein2011} for collections of sets.

\section{Conclusions}
\label{sec::conclusions}

In this survey, we have discussed polyhedral computing and polyhedral control design. We argue that certain classes of polytopes, specifically those with low complexity and balanced facet-to-vertex ratios, offer unique advantages for modeling and control. These geometric constructs, when combined with convex optimization formulations, yield scalable and powerful tools for tackling a wide range of set-based control tasks.

The historical origins of polyhedral computing, including the seminal contributions of Minkowski, Gale, Weyl, and Sommerville, as well as McMullen's groundbreaking polytope algebra and Motzkin's double description method, date back to the dawn of the previous century. This rich heritage stands in stark contrast to the more recent emergence of polyhedral control design as a distinct and vibrant field. While pioneering efforts by scholars like Bitsoris, Gutman, Cwickel, and Blanchini in the 1970s, 80s, and 90s laid the foundational stones, it is the explosive growth of modern convex optimization tools over the past two decades that has fundamentally reshaped the landscape.

This transformation has been fueled by innovative applications of convex linear programming, exemplified by Trodden's approach to computing positive invariant polytopes and Kerrigan's set-theoretic classification of convex robust control tubes. Additionally, the field has witnessed a proliferation of novel convex Tube MPC formulations, thanks to the seminal works of Chisci, Goulart, Kerrigan, Langson, Mayne, Raković, and an impressive list of other researchers.

As we look into the future, it is evident that the synergies between polyhedral computing and convex optimization will persist as a catalyst for groundbreaking progress in set-based control. This synergy promises to deepen our comprehension of the relationship between structural intricacies and computational complexity within set-based control frameworks, thereby facilitating the application of these tools to tackle increasingly complex, higher-dimensional challenges.

\section*{Acknowledgments}

The authors are grateful to Mario Zanon for his comments on this article.

\bibliographystyle{plain}
\bibliography{references.bib}

\end{document}